\documentclass[12pt,british]{amsart}

\usepackage[latin1]{inputenc}
\usepackage[T1]{fontenc}
\usepackage{babel,eucal,url,amssymb,stmaryrd,booktabs,%
enumerate,amscd,paralist}
\usepackage{mathrsfs}

\theoremstyle{plain}
\newtheorem{theorem}{Theorem}[section]
\newtheorem{proposition}[theorem]{Proposition}
\newtheorem{lemma}[theorem]{Lemma}
\newtheorem{corollary}[theorem]{Corollary}

\theoremstyle{definition}
\newtheorem{definition}[theorem]{Definition}

\theoremstyle{remark}
\newtheorem{remark}[theorem]{Remark}

\newtheorem*{acknowledgements}{Acknowledgements}

\numberwithin{equation}{section}

\DeclareMathOperator{\alg}{alg}
\DeclareMathOperator{\DT}{DT}

\DeclareMathOperator{\dist}{dist}
\DeclareMathOperator{\ex}{e}
\DeclareMathOperator{\im}{i}

\DeclareMathOperator{\F2}{\mathbf F_2}

\newcommand{\norm}[1]{\left\lVert #1\right\rVert}

\begin{document}

\title[Free entropy dimension of DT-operators]{The
non-microstates free entropy \\ dimension of DT-operators} 
 
\author{Lars Aagaard}
\address[Lars Aagaard]{Department of Mathematics and Computer Science\\
University of Southern Denmark\\
Campusvej 55\\
DK-5230 Odense M\\
Denmark}
\email{laa@imada.sdu.dk}


\begin{abstract}
Dykema and Haagerup introduced the class of DT-operators \cite{DT}
and also showed that every DT-operator  
generate $L(\F2)$ \cite{invsub}, the von 
Neumann algebra generated 
by the free group on two generators. 
In this paper we prove that Voiculescu's non-microstates free entropy
dimension is 2 for all DT-operators.   
\end{abstract}


\maketitle


\section{Introduction}

The class of $\DT$-operators was introduced by Dykema and Haagerup in
 \cite{DT}.  
 For $n\in
\mathbb N$ let $\mathcal M_n(\mathbb C)$ be the space of $n\times
n$-matrices with entries being random varibles over a classical
 probability space having moments of all
orders. Let $\tau_n:\mathcal M_n(\mathbb C)\to \mathbb C$ be the
expectation of the normalized trace on $\mathcal M_n(\mathbb 
C)$. Then $(\mathcal M_n(\mathbb C),\tau_n)$ is a $*$-noncommutative
probability space.

Let $\mu$ be a compactly supported Borel-measure on $\mathbb C$ and
let $c>0$. Let $D_n\in  
\mathcal M_n(\mathbb C)$ be diagonal matrices with i.i.d diagonal entries
all having distribution $\mu$. Let $T_n\in
\mathcal M_n(\mathbb C)$ be strictly 
upper triangular matrices such that the $n(n-1)$ real and imaginary
parts of the entries
above the diagonal consists of a family of i.i.d. centered gaussian
random variables with variance $1/{2n}$. Define $Z_n=D_n+cT_n$.
    Let $Z$ be an element in a $*$-non-commutative
probability space, $(\mathcal M,\tau)$. $Z$ is a
$\DT(\mu,c)$-element if its $*$-moments, 
$\tau(Z^{\epsilon_1} Z^{\epsilon_2}\cdots Z^{\epsilon_k})$, is
determined by the the limit
\begin{equation} \label{DT-limit}
\lim_{n\to \infty}\tau_n(Z_n^{\epsilon_1} Z_n^{\epsilon_2}\cdots
Z_n^{\epsilon_k})
\end{equation}
for all $k\in \mathbb N$ and $\epsilon_1,\ldots,\epsilon_k\in
\{1,*\}$. The limit (\ref{DT-limit}) always exists \cite[Th. 2.1]{DT}.

A  $\DT(\mu,c)$-operator is a $\DT(\mu,c)$-element constructed in a
$W^*$-probability space, and a $\DT$-operator is a 
$\DT(\mu,c)$-operator for some $\mu$ and $c$. Haagerup and Dykema has shown
  \cite{invsub} that actually every $\DT$-operator generate a von Neumann
  algebra isomorphic to $L(\F2)$, the
von Neumann algebra generated by the free group on two generators.

Voiculescu has introduced two kinds of entropy; a microstates free entropy,
$\chi$, and 
a non-microstates free entropy, $\chi^*$, and he conjectures that
$\chi=\chi^*$. Following the definition of the microstates free entropy
dimension, $\delta$, from $\chi$ we define the non-microstates free
entropy dimension, $\delta^*$, similarly from $\chi^*$. So if $\chi=\chi^*$
then of course $\delta=\delta^*$. 

\'Sniady has shown a formula \cite{SniadyDT} for the microstates free
entropy of $\DT$-operators. In particular it follows from his results that a 
$\DT(\delta_0,1)$-operator or upper triangular operator which is also just 
the limit in $*$-moments of the $T_n's$ above, has microstates 
free entropy $-\infty$. This makes it an interesting problem to compute the free entropy dimension of this operator since $-\infty$-microstates entropy is the only chance that the microstates entropy dimension can be less than 2. We prove the following theorem
\begin{theorem}
Let $\mu$ be a compactly supported Borel measure on $\mathbb C$, let
$c>0$ and let $Z$ be a $\DT(\mu,c)$-operator. Then
\begin{equation}
\label{eq:dimZequal2}
\delta^*(Z)=2.
\end{equation}
\end{theorem}

We also consider $\delta^*(\cdot,\mathcal B)$, the non-microstates
free entropy dimension with respect to an algebra, $\mathcal B$, 
and we show that $\delta^*(a_1,\ldots,a_n:\mathcal B)$ of self-adjoint
variables $a_1,\ldots,a_n$ can only be different from $n$ if
the non-microstates free Fischer information,
$\Phi^*(a_1,\ldots,a_n:\mathcal B)$, is $+\infty$.


\begin{acknowledgements}
It is a great pleasure for me to thank my advisor Professor Uffe
Haagerup for suggesting the problem to me and for stimulating and
fruitful discussions on the subject. 
I would also like to thank Peter Ainsworth for letting me publish the
results in the appendix in this paper. Finally I would like to thank
A. Nica for pointing out the reference \cite{NSS2}. 
\end{acknowledgements}

\section{Review on free Fischer information}

The following review on free Fisher information can be read out of
\cite{NSS1}, \cite{NSS2} and \cite{Combinatorics}.  Free Fischer
information was originally introduced by Voiculescu in \cite{Voi1},
\cite{Voi5} and further investigated in \cite{Voi6}. 

In this section we will, unless otherwise stated, let $(\mathcal
A,\phi)$ be a
$W^*$-probability space with $\phi$ a 
faithful, normal trace, let $1\in \mathcal B\subset A$ be a unital
$W^*$-sub-algebra. Also we let $E_\mathcal B:\mathcal A\to \mathcal B$ be a
conditional expectation such that $(\mathcal B\subset \mathcal
A,E_\mathcal B)$ is a $\mathcal B$-probability space compatible to
$(\mathcal A,\phi)$ in the sense that $\phi=\phi\circ E_\mathcal B$. Let
$L^2(\mathcal A,\phi)$ be the Hilbert space 
completion of $\mathcal A$ with respect to the norm $\norm{a}_\phi =
\sqrt{\phi(a^*a)}$ for $a\in \mathcal A$.

 If  $\mathcal X$ is a subset of
$\mathcal A$ then $L^2(\mathcal X,\phi)$ will denote the $L^2$-completion,
$\overline{\textrm{alg}(\mathcal X,\mathcal X^*)}^{||\cdot 
||_\phi}$, of the unital $*$-algebra generated by $\mathcal X$. 

We define a self-adjoint family, $(a_i)_{i\in I} \subset \mathcal
A$, to be a family of operators such that for all $i\in I$ there exists
$j\in I$ such that $a_i^*=a_j$.

\begin{definition}
Let $(\mathcal A,\phi)$ be a $W^*$-probability space with $\phi$ a
faithful trace. Let $1\in \mathcal B \subset A$ be a unital
$W^*$-sub-algebra, and let $(a_i)_{i\in I}$ be a self-adjoint family of
random variables in $\mathcal A$. Then a family of vectors $(\xi)_{i\in
I}$ from $L^2(\mathcal A,\phi)$ fulfills the conjugate relations for
$(a_i)_{i\in I}$ with respect to $\mathcal B$ if
\begin{multline}
\label{conjugate}
\phi(\xi_i b_0a_{i_1}b_1 a_{i_2}\cdots a_{i_n}b_n) \\ = \sum_{m=1}^n
\delta_{i,i_m}\phi(b_0a_{i_1}\dots a_{i_{m-1}}b_{m-1})\phi(b_m
a_{i_{m+1}}b_{m+1} \cdots a_{i_n}b_n),
\end{multline} 
for every $n\geq 0$, $b_0,b_1,\ldots,b_n\in B$ and $i,i_1,i_2,\ldots
i_n\in I$.

A family of vectors $(\xi_i)_{i\in I}\subset L^2(\mathcal A,\phi)$ is
said to be a conjugate  
system for a self-adjoint family of operators, $(a_i)_{i\in I}\subset
\mathcal A$, with
respect to $\mathcal B$ if it satisfies 
the conjugate relations (\ref{conjugate}) and if furthermore $(\xi_i)_{i\in I}
\subset L^2((a_i)_{i\in I}\cup \mathcal B,\phi)$.  
\end{definition}

\begin{remark} \label{remark1.2}
\begin{enumerate}
\item[(a)] The above definition is to be understood as $\phi(\xi_ib) = 0$ for
all $i\in I$. Since $\mathcal B$ is unital we thus have $\phi(\xi_i)=0$
for all $i\in I$. 
\item[(b)] If a conjugate system $(\xi_i)_{i\in I}$ exists then it is
unique since (\ref{conjugate}) is a prescription for taking inner
products with monomials of the form $b_0a_{i_1}b_1a_{i_2}\cdots a_{i_n}b_n$ from $L^2((a_i)_{i\in I}\cup \mathcal B,\phi)$, so the inner product of an element from $(\xi_i)_{i\in I}$ with an arbitrary element from $L^2((a_i)_{i\in I}\cup \mathcal B,\phi)$ is completely determined.
\item[(c)] If one can find $(\xi)_{i \in I}$ that fulfills the conjugate
relations, (\ref{conjugate}), for a self-adjoint family, $(a_i)_{i\in
I}\subset \mathcal A$, with respect to a unital sub-algebra $\mathcal B$
of $\mathcal A$, then if $P:L^2(\mathcal A,\phi)\to L^2((a_i)_{i\in
I}\cup \mathcal B,\phi)$ is the Hilbert space projection then $(P\xi_i)_{i\in I}$ is a conjugate system for $(a_i)_{i\in I}$ with respect to $\mathcal B$.  
\item[(d)] If $(a_i)^n_{i=1}$ are all self-adjoint then Voiculescu
originally denoted the conjugate variables of $(a_i)^n_{i=1}$ with
respect to $\mathcal B$ by  
\begin{equation*}
\mathscr J(a_i:\mathcal B[a_1,\ldots,a_{i-1},a_{i+1},\ldots, a_n])
\end{equation*}
 for $i\in \{1,\ldots,n\}$. 
\end{enumerate}
\end{remark}

\begin{definition}
  Let $(\mathcal A,\phi)$ be a $W^*$-probability space with $\phi$ a
faithful trace. Let $\mathcal B \subset A$ be a unital
$W^*$-sub-algebra, and let $(a_i)_{i\in I}$ be a self-adjoint family of
random variables in $\mathcal A$. If $(a_i)_{i\in I}$ has a conjugate system, $(\xi)_{i\in
I}$, with respect to $\mathcal B$ then we define the free Fischer information of $(a_i)_{i\in I}$ with respect to $B$ as:
\begin{equation}
  \label{fischer}
  \Phi^*((a_i)_{i\in I}:\mathcal B) = \sum_{i\in I} ||\xi_i||_\phi^2.
\end{equation}
If no conjugate system exists for $(a_i)_{i\in I}$ with respect to $\mathcal B$ we define $\Phi^*((a_i)_{i\in I}:\mathcal B) = +\infty$. If $\mathcal B = \mathbb C 1$ then we define 
$$
\Phi^*((a_i)_{i\in I}) = \Phi^*((a_i)_{i\in I}: \mathbb C 1),
$$
and we call this the free Fischer information of $(a_i)_{i\in I}$.
\end{definition}

\begin{remark} \label{fischerremark}
  
\begin{enumerate}
\item[(a)] If $(\xi_i)_{i\in I}$ satisfies the conjugate relations for a
self-adjoint family $(a_i)_{i\in I}\subset \mathcal A$ with respect to
$\mathcal B$ and if $P:L^2(\mathcal A,\phi)\to L^2((a_i)_{i\in I}\cup
\mathcal B,\phi)$ is the Hilbert space projection then by remark
\ref{remark1.2} (c) we know that $(P\xi_i)_{i\in I}$ is a conjugate
system for $(a_i)_{i\in I}$ with respect to $\mathcal B$ and since
projections are norm decreasing we conclude that 
  \begin{equation}
    \label{relations}
    \Phi^*((a_i)_{i\in I}:\mathcal B) \leq \sum_{i\in I} ||\xi_i||^2_\phi.
  \end{equation}
\item[(b)] If $r\in \mathbb R$ is a strictly positive scalar then one
easily sees that if $(a_i)_{i\in I}\subset \mathcal A$ is a self-adjoint
family of random variables with conjugate system $(\xi_i)_{i\in I}$ then $(\frac{1}{r}\xi_i)_{i\in I}$ is a conjugate system for $(ra_i)_{i\in I}$ and thus
\begin{equation}
\label{scalarfischer}
\Phi^*((ra_i)_{i\in I}:\mathcal B) = \frac{1}{r^2} \Phi^*((a_i)_{i\in I}:\mathcal B). 
\end{equation}
\item[(c)] The free Fischer information respects inclusion of
sub-algebras in the following sense. If $(\mathcal A,\phi)$ is a
$W^*$-probability space and $1\in \mathcal B_1\subset \mathcal B_2$ are
to unital $W^*$-sub-algebras then if $(a_i)_{i\in I}\subset \mathcal A$
is a self-adjoint system then  
 \begin{equation}
    \label{inclusion}
    \Phi^*((a_i)_{i\in I}:\mathcal B_1)\leq \Phi^*((a_i)_{i\in I}:\mathcal B_2)
  \end{equation}
because if a conjugate system for $(a_i)_{i\in I}$ exists with respect to
$\mathcal B_2$ then this conjugate system will also satisfy the
conjugate relations for $(a_i)_{i\in I}$ with respect to $\mathcal B_1$,
and hence (\ref{inclusion}) follows from (\ref{relations}) in (a). 
\end{enumerate}
\end{remark}

The following theorem is a special case of \cite[Th. 4.1]{NSS2}. 
Concerning cumulants we adopt the tensor-product notation of \cite{AMS}.

\begin{theorem}\cite[Th. 4.1]{NSS2}\label{cumulantconjugate}
  Let $(\mathcal A,\phi)$ be a $W^*$-probability space where $\phi$ is a
  faithful, normal trace and let $(\mathcal B\subset \mathcal
  A,E_\mathcal B)$ be a $\mathcal B$-probability space compatible to
  $(\mathcal A,\phi)$. Let $(a_i)_{i\in I}$ be a 
  self-adjoint family of random variables in $\mathcal A$. Then
  $(\xi_i)_{i\in I}$ satisfies the conjugate relations for $(a_i)_{i\in
  I}$ with respect to $\mathcal B$ if and only if 
  \begin{equation}
    \label{cumconj}
    \kappa^\mathcal B_{n+1}(\xi_i\otimes_\mathcal B b_0a_{i_1}\otimes_\mathcal B b_1a_{i_2}\otimes_\mathcal B\cdots \otimes_\mathcal B b_{n-1}a_{i_n})=
    \begin{cases}
      \delta_{ii_1}\phi(b_0)1 & \textrm{for } n=1 \\
      0 & \textrm{for } n\neq 1.
    \end{cases}
  \end{equation}
for all $b_0,\ldots, b_{n-1}\in \mathcal B$ and $i,i_1,\ldots i_n\in I$.
\end{theorem}

\begin{remark} \label{rem1.6}
  
Consider a non-self-adjoint random variable $a\in \mathcal A$.
Then a conjugate system for $(a,a^*)$ must have the form $(\xi,\xi^*)$
because of the tracial properties of $\phi$ and the
conjugate relations, (\ref{conjugate}).  
From theorem \ref{cumulantconjugate} it is easy to see that
$\left(\frac{\xi+\xi^*}{\sqrt{2}}, -\frac{\xi-\xi^*}{\im\sqrt{2}}\right)$ is a
conjugate system for
$\left(\frac{a+a^*}{\sqrt{2}},\frac{a-a^*}{\im\sqrt{2}}\right)$ with
respect to $\mathcal B$ so we conclude that
\begin{eqnarray*}
\label{}
\Phi^*(a,a^*) & = & ||\xi||^2_\phi + ||\xi||^2_\phi = \phi(\xi^*\xi) +
\phi(\xi\xi^*) \\
& =  & \phi\left(\tfrac{\xi^2+ (\xi^*)^2 + \xi^*\xi + \xi\xi^*}{2}\right) 
 + \phi\left(-\tfrac{\xi^2+ (\xi^*)^2 - \xi^*\xi - \xi\xi^*}{2}\right) \\
& = & \left|\left|\tfrac{\xi+\xi^*}{\sqrt{2}}\right|\right|_\phi^2 +
\left|\left|-\tfrac{\xi-\xi^*}{\im\sqrt{2}}\right|\right|_\phi^2 = \Phi^*\left(\tfrac{a+a^*}{\sqrt{2}},\tfrac{a-a^*}{\im\sqrt{2}}\right).
\end{eqnarray*} 
Combining with (\ref{scalarfischer}) of remark
\ref{fischerremark} we have 
$$2\Phi^*(a,a^*)  = \Phi^*(\Re a, \Im
a).$$

\end{remark}

\section{Non-microstates free entropy dimension} \label{dimsection}

The non-microstates free entropy for several self-adjoint random variables was
originally defined by Voiculescu.

\begin{definition} \cite[Def. 7.1]{Voi5} \label{*entropy}
Let $(a_i)_{i=1}^n\subset\mathcal A$ be a collection of self-adjoint random
variables in a $W^*$-probability space $(A,\phi)$ where $\phi$ is a
faithful normal tracial state, and let $1\in \mathcal B\subset \mathcal A$ be a
unital $W^*$-sub-algebra of $\mathcal A$. The non-microstates free
entropy of $(a_i)_{i=1}^n$ with respect to $\mathcal B$ is then defined as
\begin{multline}
\label{entropydef}
\chi^*(a_1,\ldots,a_n:\mathcal B) \\ = \frac 12 \int_0^\infty \frac{n}{1+t} -
\Phi^*(a_1+\sqrt t S_1,\ldots a_n+\sqrt t S_n)\text{d} t  + \frac n 2 \log(2\pi \ex),
\end{multline}
where $S_1,\ldots,S_n$ are standard semicircular elements such that $\{S_1\},\ldots,\{S_n\}$ and $\{\alg((a_i)_{i=1}^n\cup \mathcal B)\}$
are free sets.
\end{definition}

The following property of $\chi^*$ is shown by Voiculescu in \cite{Voi5}.


\begin{proposition}\cite[Prop 7.2]{Voi5} \label{upperchibound}
Let $(a_i)_{i=1}^n, \mathcal B,\mathcal
A$, $(S_i)_{i=1}^n$ and $\phi$ be as in definition \ref{entropydef}.
Let $C^2= \phi(a_1^2+\cdots a_n^2)$. Then
\begin{equation} \label{chiovre}
\chi^*(a_1,\ldots,a_n:\mathcal B) \leq \tfrac{n}{2}\log(2\pi \ex n^{-1}C^2).
\end{equation}
\end{proposition}

In \cite{Voi2} Voiculescu defined the microstates free entropy dimension. 
We consider the non-microstates free analogue.
\begin{definition} \label{dimension}
Let $(a_i)_{i=1}^n\subset \mathcal A$ be a collection of self-adjoint random
variables in a $W^*$-probability space $(A,\phi)$ where $\phi$ is a
faithful normal tracial state, and let $1\in \mathcal B\subset \mathcal A$ be a
unital $W^*$-sub-algebra of $\mathcal A$. The non-microstates free entropy
dimension of 
$(a_i)_{i=1}^n$ with respect to $\mathcal B$ is defined by
\begin{equation}
\label{dimensioneq}
\delta^*(a_1,\ldots,a_n:\mathcal B) = n+ \limsup_{\epsilon\to 0}
\frac{\chi^*(a_1+\epsilon S_1,\ldots,a_n+\epsilon S_n:\mathcal B)}{|\log\epsilon|},
\end{equation}
where $S_1,\ldots,S_n$ are standard semicircular variables such that 
$\{S_1\},\ldots,\{S_n\}$ and $\{(a_i)_{i=1}^n\cup\mathcal B\}$ are $\phi$-free. 
\end{definition} 

An easy upper bound of the non-microstates free entropy dimension
follows from proposition \ref{upperchibound}
\begin{proposition} \label{deltaovre}
Let $(a_i)_{i=1}^n\subset \mathcal A$ be a collection of not all zero
self-adjoint random 
variables in a $W^*$-probability space $(A,\phi)$ where $\phi$ is a
faithful normal tracial state, and let $1\in \mathcal B\subset \mathcal A$ be a
unital $W^*$-sub-algebra of $\mathcal A$. Let $(S_i)_{i=1}^n$ be a standard semicircular family $\phi$-free from $(a_i)_{i=1}^n$. Then
\begin{equation}
\delta^*(a_1,\ldots,a_n:\mathcal B)\leq n.
\end{equation}
\end{proposition}
\begin{proof}
Let $C^2= \phi(a_1^2+\cdots a_n^2)>0$. Then 
\begin{equation*}
\phi((a_1+\epsilon S_1)^2+\cdots +(a_n+\epsilon S_n)^2) =
C^2+n\epsilon^2, 
\end{equation*}
so we infer from (\ref{chiovre}) of proposition \ref{upperchibound}
that 
\begin{eqnarray*}
\delta^*(a_1,\ldots,a_n:\mathcal B) &  = & n + \limsup_{\epsilon\to 0}
\frac{\chi^*(a_1+\epsilon S_1,\ldots,a_n+\epsilon S_n:\mathcal
B)}{|\log\epsilon|} \\ & \leq & n + \limsup_{\epsilon\to
 0^+}\frac{\tfrac{n}{2}\log(2\pi \ex n^{-1} (C^2+\epsilon^2 n))}{|\log
\epsilon|} \\ & = & n,
\end{eqnarray*}
since $C^2$ is strictly positive.
\end{proof}

A lower bound of the non-microstates free entropy dimension can be
estimated from the non-microstates free Fischer information in the following 
sense.
\begin{proposition} \label{lowerdeltabound}
Let $(a_i)_{i=1}^n\subset \mathcal A$ be a collection of self-adjoint random
variables in a $W^*$-probability space $(A,\phi)$ where $\phi$ is a
faithful normal tracial state, and let $1\in \mathcal B\subset \mathcal A$ be a
unital $W^*$-sub-algebra of $\mathcal A$. The non-microstates 
free entropy dimension of
$(a_i)_{i=1}^n$ with respect to $\mathcal B$ is bounded from below by
\begin{equation}
\label{lowerdimbound}
\delta^*(a_1,\ldots,a_n:\mathcal B) \geq n - \limsup_{t\to 0^+}
(t\Phi^*(a_1+\sqrt{t}S_1,\ldots,a_n+\sqrt{t}S_n:\mathcal B)),
\end{equation} 
where $S_1,\ldots,S_n$ is a family of standard semicircular random
variables such that 
$\{S_1\},\ldots,\{S_n\}$ and $\{\mathcal B\cup(a_i)_{i=1}^n\}$ are
$\phi$-free. If the $\limsup$ on the right-hand side of
(\ref{lowerdimbound}) is convergent then inequality in
(\ref{lowerdimbound}) becomes equality.
\end{proposition}

\begin{proof} We only show the inequality of (\ref{lowerdimbound}). The
last statement about equality follows by trivial modifications of the
argument.
Let $\mathcal A,\mathcal B,\phi, (a_i)_{i=1}^n$ and $(S_i)_{i=1}^n$ be
as in the proposition. Define
\begin{equation*}
a := \limsup_{t\to
0^+}(t\Phi^*(a_1+\sqrt{t}S_1,\ldots,a_n+\sqrt{t}S_n:\mathcal B)),
\end{equation*}
and let $\epsilon > 0$. There exists an interval $\delta>0$ such that
\begin{equation*}
\Phi^*(a_1+\sqrt{t}S_1,\ldots,\sqrt{t}S_n:\mathcal B) \leq \frac{a+\epsilon}{t},
\end{equation*}
for all $0<t<\delta$.

Using definition \ref{*entropy} carefully implies
\begin{multline}
\chi^*(a_1+\sqrt{\delta}S_1,\ldots,a_n+\sqrt{\delta}S_n:\mathcal B)
\\ -\chi^*(a_1+\sqrt{t}S_1,\ldots,a_n+\sqrt{t}S_n:\mathcal B) \\
= \frac 1 2 \int_t^\delta \frac{n}{1+s-\delta} - \Phi^*
(a_1+\sqrt{s}S_1,\ldots,a_n+\sqrt{s}S_n:\mathcal B)\text{d} s \\ 
+\frac 1 2\int_\delta^\infty 
\frac{n}{1+s-\delta} - \frac{n}{1+s-t}\text{d} s\\
= \frac 1 2 \int_t^\delta
 \Phi^*(a_1+\sqrt{s}S_1,\ldots,a_n+\sqrt{s}S_n:\mathcal B)\text{d} s\\
 \leq \frac 1 2 \int_t^\delta \frac{a+\epsilon}{s}\text{d} s =
 \frac{a+\epsilon}{2}\log\left(\frac{\delta}{t}\right), 
\end{multline}
for all $0<t<\min\{\delta,1\}$, so
\begin{multline*}
\chi^*(a_1+\sqrt{t}S_1,\ldots,a_n+\sqrt{t}S_n:\mathcal B) \\
\geq \chi^*(a_1+\sqrt{\delta}S_1,\ldots,a_n+\sqrt{\delta}S_n:\mathcal B)
- \frac{a+\epsilon}{2}\log\delta - \frac{a+\epsilon}{2}|\log t |.
\end{multline*}
From this we deduce that
\begin{equation*}
\liminf_{t\to 0^+}
\frac{\chi^*(a_1+\sqrt{t}S_1,\ldots,a_n+\sqrt{t}S_n:\mathcal B)}{|\log
t|} \geq -\frac{a+\epsilon}{2},
\end{equation*}
and since $\epsilon$ was chosen arbitrarily comparing to
(\ref{dimensioneq}) implies that 
\begin{eqnarray*}
\delta^*(a_1,\ldots,a_n:\mathcal B) &= & n + 2\limsup_{t\to 0^+}
\frac{\chi^*(a_1+\sqrt{t}S_1,\ldots,a_n+\sqrt{t}S_n:\mathcal B)}{|\log
t|} \\
& \geq & n +2\liminf_{t\to 0^+}\frac{\chi^*(a_1+\sqrt{t}S_1,\ldots,a_n+\sqrt{t}S_n:\mathcal B)}{|\log
t|} \\ 
& \geq & n - a \\
& =&  n- \limsup_{t\to
0^+}(t\Phi^*(a_1+\sqrt{t}S_1,\ldots,a_n+\sqrt{t}S_n:\mathcal B)).
\end{eqnarray*}
\end{proof}

We remark that it might be tempting to adopt the right-hand side of
(\ref{lowerdimbound}) as the 
definition of the non-microstates free entropy dimension with respect to
$\mathcal B$. Voiculescu has shown the following Free Stam inequality.

\begin{proposition}\cite[Prop. 6.5]{Voi5}  \label{boundPhi2}
Let $(a^1_i)_{i=1}^n,(a^2_i)_{i=1}^n \subset \mathcal A$ be a collection
of self-adjoint random
variables in a $W^*$-probability space $(A,\phi)$ where $\phi$ is a
faithful normal tracial state, and let $1\in \mathcal B_1,\mathcal
B_2\subset \mathcal A$ be a unital $W^*$-sub-algebras of $\mathcal A$. If
$\alg((a_i^1)_{i=1}^n,\mathcal B_1)$ and $\alg((a_i^2)_{i=1}^n,\mathcal
B_2)$ are $\phi$-free then
\begin{multline} \label{Philess2}
\Phi^*(a_1^1+a_1^2,\ldots,a_n^1+a_n^2:W^*(\mathcal B_1,\mathcal
B_2))^{-1} \\ \geq \Phi^*(a_1^1,\ldots,a_n^1:\mathcal B_1)^{-1} +\Phi^*(a_1^2,\ldots,a_n^2:\mathcal B_2)^{-1}
\end{multline}  
\end{proposition}


Combining the lower bound of $\delta^*$ from proposition
\ref{lowerdeltabound}, proposition \ref{deltaovre}  and proposition
\ref{boundPhi2} we have the following interesting which limits
the free entropy dimension problem to the case of infinite free Fischer
information. 
\begin{corollary} \label{infinitePhi}
  Let $(a_i)_{i=1}^n\subset \mathcal A$ be a collection of self-adjoint random
variables in a $W^*$-probability space $(A,\phi)$ where $\phi$ is a
faithful normal tracial state, and let $1\in \mathcal B\subset \mathcal A$ be a
unital $W^*$-sub-algebra of $\mathcal A$. Then
\begin{equation}
\Phi^*(a_1,\ldots,a_n:\mathcal B) < \infty \Rightarrow \delta^*(a_1,\ldots,a_n:\mathcal B) = n
\end{equation}
\end{corollary}
\begin{proof}
  Apply the Free Stam inequality (\ref{Philess2}) on the algebras
  $\alg((a_i)_{i=1}^n,\mathcal B)$ and $\alg((S_i)_{i=1}^n,\mathbb C)$ and observe that
  \begin{equation*}
    \limsup_{t\to 0^+}  t\Phi^*(a_1+\sqrt{t}S_1,\ldots, a_n+\sqrt{t}S_n:\mathcal B) \leq \limsup_{t\to 0^+}\frac{nt}{\tfrac{n}{\alpha}+t} = 0,
  \end{equation*}
when $\alpha= \Phi^*(a_1,\ldots,a_n:\mathcal B)< \infty$.
\end{proof}

\section{$\chi^*$ and $\delta^*$ for a non-self-adjoint random variable}

We translate the results from section \ref{dimsection} into the context of one 
non-self-adjoint variable by considering real and imaginary parts. Let
$(\mathcal A,\phi)$ be $W^*$-probability space with a normal 
faithful tracial state, and let $a\in \mathcal A$ be non-self-adjoint
random variable. 
For the microstates free entropy we have $\chi(a) = \chi(\Re a,\Im a)$
\cite[Prop. 6.5.5]{HiaiPetz} so it seems natural to define
\begin{equation} \label{entropystar}
\chi^*(a:\mathcal B) := \chi^*(\Re a, \Im a:\mathcal B).
\end{equation}

%

and thus also 
\begin{equation} \label{dimensionstar}
\delta^*(a:\mathcal B) := \delta^*(\Re a, \Im a:\mathcal B)
\end{equation}


Using proposition \ref{lowerdeltabound} we have the following corollary on
non-microstates free entropy dimension of a single non-self-adjoint
variable. 

\begin{corollary} \label{nonsavurdering}
Let $a \in\mathcal A$ be a non-self-adjoint family of random
variables in a $W^*$-probability space $(A,\phi)$ where $\phi$ is a
faithful normal tracial state, and let $1\in \mathcal B\subset \mathcal A$ be a
unital $W^*$-sub-algebra of $\mathcal A$. Then a lower bound of the
non-microstates free entropy dimension of $a$ with respect to $\mathcal
B$ is 
\begin{equation} \label{deltavurdering}
\delta^*(a:\mathcal B) \geq 2 - \limsup_{t\to
0^+}(t\Phi^*(a+\sqrt{t}Y,(a+\sqrt{t}Y)^*:\mathcal B)),
\end{equation}
where $Y$ is a standard circular random variable $*$-free from $a$ and
$\mathcal B$. If the $\limsup$ on the right-hand side of (\ref{deltavurdering}) is convergent then inequality in (\ref{deltavurdering}) becomes equality.
\end{corollary}
 
\begin{proof}
Just use (\ref{lowerdimbound}) of proposition \ref{lowerdeltabound} on the 
self-adjoint variables $\Re a$ and $\Im a$ and substitute $S_1$ and $S_2$
with a single circular random variable $Y$ $*$-free from $a$ such that
$\frac{Y+Y^*}{\sqrt 2}=S_1$ and $\tfrac{Y-Y^*}{\im\sqrt 2}=S_2$. Then
(\ref{lowerdimbound}) becomes
\begin{eqnarray*}
\delta^*(a:\mathcal B) & = & \delta^*(\Re a,\Im a:\mathcal B) \\ 
 & \geq  &  2 - \limsup_{t\to 0^+} (t\Phi^*(\Re
 (a+\sqrt{2t}Y),\Im(a+\sqrt{2t}Y):\mathcal B) \\
\end{eqnarray*}
Now remark \ref{rem1.6} implies (\ref{deltavurdering}).

\end{proof}

\section{Properties of DT-operators} \label{DT-section}

In this section we first review some basic properties of DT-operators which 
can be read out of \cite{DT} and \cite{invsub}. We then deduce some
results relating several DT-operators and standard 
circular operators.

From now on, unless otherwise stated, we let $\mathcal M$ be a von
Neumann algebra equipped with a faithful, normal, 
tracial state, $\tau$, such that $(\mathcal M,\tau)$ is a
$W^*$-probability space. We also assume that we are given an injective, unital,
normal $*$-homomorphism $\lambda: L^\infty[0,1]\to \mathcal M$, such
that $\tau\circ\lambda(f)=\int_0^1 f(t)\text{d} t$ for $f\in L^\infty[0,1]$. We let
$\mathcal D$ be the picture of $L^\infty[0,1]$ in $\mathcal M$ and let
$E_\mathcal D:\mathcal M 
\to \mathcal D$ be the trace-preserving conditional expectation onto
$\mathcal D$. We will identify $\mathcal D$ 
and $L^\infty[0,1]$ and thus consider elements of $\mathcal D$ as
$L^\infty$-functions. Assume also that $(X_i)_{i=1}^\infty\subset
\mathcal M$ is a 
standard semicircular family of random variables $\tau$-free from 
$\mathcal D$. Then $(W^*(\mathcal D\cup
\{X_i\}))_{i=1}^\infty$ are $\mathcal D$-free with respect to
$E_\mathcal D$. As in
~\cite{DT} (we also adopt the notation of ~\cite{DT}) we can now construct
upper triangular operators, $T_i=\mathcal{UT}(X_i,\lambda)$ for
$i\in \mathbb N$. If $D_0:x\mapsto x\in \mathcal D$ then
$T_i=\mathcal{UT}(X_i,\lambda)$ is the norm limit of 
\begin{equation} \label{normDTlimit}
T_n^{(i)} = \sum_{j=1}^n 1_{\left[\tfrac{j-1}{n},\tfrac{j}{n}\right]}(D_0)X_i
1_{\left[\tfrac{j}{n},1\right]}(D_0)
\end{equation}
as $n\to \infty$ \cite[Lemma 4.1]{DT}.
It follows that $T_i\in W^*(\mathcal D\cup\{X_i\})$ for all $i\in \mathbb
N$ so $(T_i)_{i=1}^\infty$ is a $\mathcal D$-free family.
 $(T_i,T_i^*)$ is a centered $\mathcal D$-Gaussian
pair for all $i\in \mathbb N$ \cite[Appendix A]{invsub} 
and it follows that $(T_i,T_i^*)_{i=1}^\infty$ is a centered $\mathcal
D$-Gaussian set in the sense of Speicher \cite[Def. 4.2.3]{AMS}. 

Define for $f\in \mathcal D$
\begin{eqnarray}
L^*(f):x\mapsto \int_0^x f(t)\text{d} t & \textrm{and} & L(f):x \mapsto
\int_x^1 f(t)\text{d} t
\end{eqnarray}
From the appendix of
~\cite{invsub} it follows that the 
covariances are given by 
\begin{lemma}\cite[Appendix]{invsub}
Let $f\in \mathcal D$. Then
\begin{eqnarray}
E_\mathcal D (T_i f T_i^*) =  L(f) & \textrm{and} & E_\mathcal D (T_i^*
f T_i)= L^*(f), \\
 E_\mathcal D (T_i f T_i) =  0 & \textrm{and} & E_\mathcal D (T_i^*
f T_i^*)= 0,
\end{eqnarray}
and $(T_i)_{i=1}^\infty$ is a $\mathcal D$-$*$-free family with respect to
$E_\mathcal D$.
\end{lemma}

A $\DT(\mu,c)$-operator, as defined in the introduction, where $\mu$ is a
compactly supported Borel measure on $\mathbb C$ and $c>0$, can be
realized in the $W^*$-probability space $(\mathcal M,\tau)$. 
\begin{theorem}\cite[Th. 4.4]{DT} \label{theorem4.4}
Let $T_1=\mathcal{UT}(X_1,\lambda)$ where $X_1$ and $\lambda$ is as above,
and let $c>0$. Let $f\in L^\infty[0,1]$ and define $D=\lambda(f)$. Then
$D+cT_1$ is a $\DT(\mu,c)$-element, where $\mu$ is the push-forward
measure of the Lebesque-measure by f.
\end{theorem}
There is a certain freedom in choosing $D\in\mathcal D$.
\begin{lemma}\cite[Lemma 6.2]{DT} \label{lemma6.2}
Let $\mu$ be a compactly supported Borel measure on $\mathbb C$. Then
there is $f\in L^\infty[0,1]$ whose distribution is $\mu$ and such that
if $D=\lambda(f)$ and if $T_1=\mathcal{UT}(X_1,\lambda)$ is as above,
then for any $c>0$ $D$ itself lies in the $W^*$-algebra generated by
$D+cT_1$. 
\end{lemma}
Dykema and Haagerup also showed that surprisingly we have 
\begin{theorem}\cite[Th. 2.2]{invsub} \label{LF2}
Let $T_1=\mathcal{UT}(X_1,\lambda)$ be as above. If $D_0:x\mapsto x\in
\mathcal D$ is the identity then $D_0\in W^*(T_1)$. 
\end{theorem}

Summarizing on theorem \ref{theorem4.4}, lemma \ref{lemma6.2} and
theorem \ref{LF2} then a
$\DT(\mu,c)$-element can be constructed as
$Z:=D+cT_1$ where the distribution of $D\in \mathcal D$ is
$\mu$, $c>0$ and $T_1=\mathcal{UT}(X_1,\lambda)$. $D$ can furthermore be 
chosen in such a way that $D\in W^*(Z)=W^*(T_1)\cong L(\F2)$.

From \cite{DT} we know that $T_i+T_i^*=X_i$ for all $i\in \mathbb N$. We
will need a relation between upper triangular operators and
standard circular variables, and for this we will use the following
theorem which is a special case of a theorem of Nica, Shlyakhtenko
and Speicher \cite{NSS3}. The
theorem as it is stated here is from \cite[Th. 6]{SnSp}.

\begin{theorem}\cite{NSS3} \label{SnSptheorem}
Let $(\mathcal M,\tau)$ be a probability space and $(\mathcal D\subset
\mathcal M,E_\mathcal D)$ a $\mathcal D$-probability space compatible to
$(\mathcal M,\tau)$ and let $\mathcal X\subset \mathcal M$. Then $\mathcal X$ and $\mathcal D$
are free in $(\mathcal M,\tau)$ iff for every $n\geq 1$ and
$x_1,\ldots,x_n\in \mathcal X$ there exists $c_n(x_1,\ldots,x_n)\in \mathbb C$
such that for every $d_1,\ldots d_{n-1}\in \mathcal D$ we have
\begin{multline*}
\kappa_n^\mathcal D(x_1d_1\otimes_\mathcal D\cdots \otimes_\mathcal D
x_{n-1}d_{n-1}\otimes_\mathcal D x_n) \\
 = c_n(x_1,\ldots,x_n)\tau(d_1)\cdots \tau(d_{n-1})1
\end{multline*}
If the above holds, we have
\begin{equation*}
c_n(x_1,\ldots,x_n)= \kappa_n^\mathbb C(x_1,\ldots,x_n).
\end{equation*}  
\end{theorem}  

\begin{lemma} \label{distributiontheorem}
   Let $T_1=\mathcal{UT}(X_1,\lambda)$ and $T_2= \mathcal{UT}(X_2,\lambda)$ 
 be upper triangular operators, and let $Y\in \mathcal M$ be a
  standard circular random variable $*$-free from $T_1$ with respect to 
$\tau$. Then for $a,b\geq 0$
  \begin{equation}
    \label{distribution}
    \sqrt{a}T_1+\sqrt{b}Y \stackrel{*-dist.}{\sim}
    \sqrt{a+b}T_1+\sqrt{b}T_2^*. 
  \end{equation} 
In particular $T_1+T_2^*$ is a circular element.
\end{lemma}
\begin{proof}
  First we prove that if $T_3=\mathcal{UT}(X_3,\lambda)$ is upper 
  triangular
  then $T_3+T_2^*$ is a standard circular variable with respect to $\tau$. We
  know that $T_3+T_2^*$ is $\mathcal D$-Gaussian since $T_2$ and $T_3$
  are $\mathcal D$-free and the $\kappa^{\mathcal D}_2$-cumulants are given by
  \begin{gather*}
  \kappa_2^\mathcal D((T_3+T_2^*)\otimes_\mathcal D d(T_3+T_2^*)^*) =  L(d)+L^*(d)
  = \tau(d)1 \\
\kappa_2^\mathcal D((T_3+T_2^*)^*\otimes_\mathcal D d(T_3+T_2^*)) =  L^*(d)+L(d)
  = \tau(d)1 \\
  \kappa_2^\mathcal D((T_3+T_2^*)\otimes_\mathcal D d(T_3+T_2^*)) =  0
  \\
 \kappa_2^\mathcal D((T_3+T_2^*)^*\otimes_\mathcal D d(T_3+T_2^*)^*) =  0
  \end{gather*}
for all $d\in \mathcal D$.
We thus see from theorem \ref{SnSptheorem} that actually $T_3+T_2^*$ is
$\mathbb C$-Gaussian with $\kappa_2^{\mathbb C}$-cumulants given by
\begin{gather*}
\kappa_2^\mathbb C(T_3+T_2^*,(T_3+T_2^*)^*) =  \kappa_2^\mathbb
C((T_2+T_3^*)^*,T_2+T_3^*) = 1  \\
\kappa_2^\mathbb C(T_3+T_2^*,T_3+T_2^*) =  \kappa_2^\mathbb
C((T_3+T_2^*)^*,(T_3+T_2^*)^*) = 0.
\end{gather*}
Hence $T_3+T_2^*$ is a standard circular random variable.
Furthermore theorem
\ref{SnSptheorem} implies that $T_3+T_2^*$ is $\tau$-free from
$\mathcal D$. Since also $T_3+T_2^*$ is $\tau$-free from $X_1$ we conclude
that $T_3+T_2^*$ is $\tau$-free from $W^*(X_1,\mathcal D)$ and in particular $T_3+T_2^*$ is $\tau$-free from $T_1$. Hence
\begin{equation*}
\sqrt{a} T_1 + \sqrt{b}(T_3+T_2^*) \stackrel{distr}{\sim} \sqrt{a}T_1+\sqrt{b}Y,
\end{equation*}
where $Y$ is a standard circular random variable $\tau$-free from
$T_1$. 

Now we just have to observe that $\sqrt{a}T_1+\sqrt{b}T_3$ is also
$\mathcal D$-Gaussian because of $\mathcal D$-freeness of $T_1$ and
$T_3$. Since 
\begin{gather*}
 \kappa_2^\mathcal D((\sqrt{a}T_1+\sqrt{b}T_3)\otimes_\mathcal D
 d(\sqrt{a}T_1+\sqrt{b}T_3)^*) =  (a+b)L(d) \\
   \kappa_2^\mathcal D((\sqrt{a}T_1+\sqrt{b}T_3)^*\otimes_\mathcal D
 d(\sqrt{a}T_1+\sqrt{b}T_3)) =  (a+b)L^*(d) \\
 \kappa_2^\mathcal D((\sqrt{a}T_1+\sqrt{b}T_3)\otimes_\mathcal D
 d(\sqrt{a}T_1+\sqrt{b}T_3)) =  0 \\
 \kappa_2^\mathcal D((\sqrt{a}T_1+\sqrt{b}T_3)^*\otimes_\mathcal D
 d(\sqrt{a}T_1+\sqrt{b}T_3)^*) =  0,
\end{gather*}
we conclude that $\sqrt{a}T_1+\sqrt{b}T_3 \stackrel{*-dist.}{\sim}
\sqrt{a+b}T_1$. Summarizing we have showed
that
\begin{equation*}
\sqrt{a}T_1+\sqrt{b}Y \stackrel{*-dist.}{\sim} \sqrt{a+b}T_1 + \sqrt{b}T_2^*,
\end{equation*}
but this is exactly what the
lemma states.
\end{proof}

The following two lemmas tells us that in the presence of the
algebra $\mathcal D$ we can ``cut out'' the lower (and upper) triangular
part of certain products and linear combinations of DT-operators.

\begin{lemma} \label{commutator} Let $D\in \mathcal D$, $c>0$. Let 
$T_1=\mathcal{UT}(X_1,\lambda)$ and $T_2=\mathcal{UT}(X_2,\lambda)$ be
upper triangular operators. Then
  \begin{equation*}
  [D+cT_1,T_2] \in \overline{W^*(T_1)([D+cT_1,T_2]+[(D+cT_1)^*,T_2^*])W^*(T_1)}^{\norm{\cdot}_\tau}.
  \end{equation*}
\end{lemma}
\begin{proof}
Let $D_0:x\to x$ be the identity-function in $\mathcal D$. By
theorem \ref{LF2} we know that $D_0\in W^*(T_1)$. Define the
projections $p_i^{(n)}
:= 1_{\left[\tfrac{i-1}{n},\tfrac{i}{n}\right]}(D_0)$ for
$i=1,\ldots,n$. If we can show that
\begin{equation} \label{combine1}
\lim_{n\to \infty}\left|\left| \sum_{\substack{i,j=1 \\i < j}}^n
p_i^{(n)}([D+cT_1,T_2]+[(D+cT_1)^*,T_2^*])p_j^{(n)} - [D+cT_1,T_2] \right|\right|_\tau=0,
\end{equation}
then the lemma follows immediately from the fact that $D_0\in
W^*(T_1)$. From the 
construction of $T_1$ in (\ref{normDTlimit}) we know that
\begin{equation} \label{sumregel}
T_i = \sum_{1\leq i\leq j\leq n} p_i^{(n)}T_ip_j^{(n)},
\end{equation}
for $i=1,2$, since we can just replace $T_i$ in (\ref{sumregel}) by 
$T_n^{(i)}$ from 
(\ref{normDTlimit})  and take the norm limit. Thus also
\begin{equation*}
T_1 T_2 = \sum_{1\leq i\leq k\leq j\leq n} p_i^{(n)} T_1 p_k^{(n)} T_2
  p_j^{(n)},
\end{equation*}
and since $p_i^{(n)}T_1 p_k^{(n)}=0$ for $i>k$ and
$p_k^{(n)}T_2p_j^{(n)}=0$ for $k>j$ we have
\begin{eqnarray*}
T_1 T_2 & = & \sum_{\substack{1\leq i\leq j\leq n \\ 1 \leq k\leq n}}
p_i^{(n)} T_1 p_k^{(n)} T_2   p_j^{(n)}. \\
 &= & \sum_{1\leq i\leq j\leq n}p_i^{(n)}T_1T_2 p_j^{(n)}.
\end{eqnarray*}

For $i \geq j$ 
\begin{equation*}
p_i^{(n)}T_1T_2p_j^{(n)} =
\sum_{k=1}^np_i^{(n)}T_1p_k^{(n)}T_2p_j^{(n)} = 
\begin{cases}
0 & i > j \\
p_i^{(n)}T_1p_i^{(n)}T_2p_i^{(n)} & i=j,
\end{cases}
\end{equation*}
because $p_i^{(n)}T_1p_k^{(n)}=0$ for $i>k$ and
$p_k^{(n)}T_2p_j^{(n)}=0$ for $j<k$.
For $i \neq j$ we have $p_i^{(n)}Dp_j^{(n)}=0$ since $D$ commutes with the orthogonal projections $p_i^{(n)}$ and $p_j^{(n)}$, so
\begin{multline*}
(D+cT_1)T_2 -\sum_{1\leq i<j\leq n} p_i^{(n)}(D+cT_1)T_2p_j^{(n)} \\ = 
\sum_{i=1}^n p_i^{(n)}(D+cT_1)p_i^{(n)}T_2p_i^{(n)} 
 =  \sum_{i,j=1}^n p_i^{(n)}(D+cT_1)p_i^{(n)}p_j^{(n)}T_2p_j^{(n)} \\
  =  \left(\sum_{i=1}^n p_i^{(n)}(D+cT_1)p_i^{(n)}\right)
\left(\sum_{j=1}^n p_j^{(n)}T_2p_j^{(n)}\right).
\end{multline*}
Reversing the role of $T_1$ and $T_2$ and subtracting we conclude that
\begin{equation} \label{combine3}
[D+cT_1,T_2] -\sum_{1\leq i< j\leq n} p_i^{(n)}[D+cT_1,T_2]p_j^{(n)} =
[D+cU^{(n)},V^{(n)}], 
\end{equation}
where $U^{(n)}=\sum_{i=1}^np_i^{(n)}T_1p_i^{(n)}$ and
$V^{(n)}=\sum_{j=1}^np_j^{(n)}T_2p_j^{(n)}$. But
$||U^{(n)}||=||V^{(n)}||=\sqrt{\tfrac{e}{n}}$ because $||T_1||=||T_2||=
\sqrt{e}$ \cite[Cor. 8.11]{DT}, so
\begin{equation} \label{combine4}
||[U^{(n)},V^{(n)}]||_\tau \leq ||[U^{(n)},V^{(n)}]|| \leq
  2\left(\norm{D}+\tfrac{\sqrt{e}}{\sqrt{n}}\right)\tfrac{\sqrt{e}}{\sqrt{n}}, 
\end{equation}

For $i<j$ we have 
\begin{equation*}
  p_i^{(n)}(D+cT_1)^*T_2^*p_j^{(n)} = \sum_{k=1}^n p_i^{(n)}(D+cT_1)^*p_k^{(n)}T_2^*p_j^{(n)} =0
\end{equation*}
because $D$ commute with the involved projections and since 
$p_i^{(n)}T_1^*p_k^{(n)}=0$ for $i<k$ and $p_k^{(n)}T_1^*p_j^{(n)}=0$ 
for $k<j$. Reversing the role of $(D+cT_1)^*$ and $T_2^*$, summing and 
subtracting we have
\begin{equation} \label{combine6}
   \sum_{1\leq i < j \leq n} p_i^{(n)}[(D+cT_1)^*,T_2^*]p_j^{(n)} =0.
\end{equation}
Hence (\ref{combine1}) follows from 
combining (\ref{combine3}), (\ref{combine4}) 
and (\ref{combine6}). 
\end{proof}

\begin{lemma} \label{cutouttheorem}
  
Let $T_1=\mathcal{UT}(X_1,\lambda)$ and $T_2=\mathcal{UT}(X_2,\lambda)$ be upper triangular operators as in the beginning of the section. Let $D\in \mathcal D$ and let $a,b\in \mathbb R_+$ be nonzero positive scalars. 
Then 
\begin{equation*}
T_1,T_2\in W^*(D+\sqrt{a}T_1+\sqrt{b}T_2^*,\mathcal D ).
\end{equation*}
\end{lemma}

\begin{proof}
We proceed as in the proof of lemma \ref{commutator},
and let $D_0:x\to x$ be the identity-function in $\mathcal D$. We only
show that $T_1\in W^*(D+\sqrt{a}T_1+\sqrt{b}T_2^*,\mathcal D )$, since
the result for $T_2$ is obtained in a similar way. Define the
projections $p_i^{(n)}
:= 1_{\left[\tfrac{i-1}{n},\tfrac{i}{n}\right]}(D_0)$ for
$i=1,\ldots,n$. Remember from the proof of lemma \ref{commutator} that
\begin{equation*}
T_1 = \sum_{1\leq i \leq j \leq n} p_i^{(n)}T_1 p_j^{(n)},
\end{equation*}
and that
\begin{equation*}
U^{(n)}:= T_1- \sum_{1\leq i < j \leq n} p_i^{(n)}T_1 p_j^{(n)} = \sum_{i=1}^n
p_i^{(n)}T_1p_i^{(n)}.
\end{equation*}
We have $||U^{(n)}||\leq \sqrt{\tfrac{e}{n}}$ and also 
\begin{equation}
  \sum_{1\leq i < j \leq n} p_i^{(n)}T_2^* p_j^{(n)} =0.
\end{equation}
We conclude that
\begin{multline*} 
\lim_{n\to \infty}\left|\left|\sqrt{a}T_1 -\sum_{\substack{i,j=1 \\i < j}}^n
p_i^{(n)}(D+\sqrt{a}T_1+\sqrt{b}T_2^*)p_j^{(n)}  \right|\right| 
\\ \leq
\lim_{n\to \infty}\sqrt{a}\norm{U^{(n)}}\leq
\lim_{n\to\infty}\sqrt{\tfrac{a\ex}{n}} 
=0.
\end{multline*} 
\end{proof}

\section{$\delta^*$ for a DT-operator with respect to $\mathcal D$} \label{deltaD}

Let $T_1=\mathcal{UT}(X_1,\lambda)$, $T_2=\mathcal{UT}(X_2,\lambda)$ be
upper triangular operators as in the beginning of section \ref{DT-section}. 
Define $Z=D+cT_1\in \mathcal M$ where $c>0$
and $D\in \mathcal D$. By theorem \ref{theorem4.4} $Z$ is a
$\DT(\mu,c)$-operator where $\mu$ is 
the distribution of $D$. Let $Y\in \mathcal M$ be
a standard circular operator $*$-free from $D$ and $T_1$ and thus also
$*$-free from $Z$.  
By corollary \ref{nonsavurdering} we are to compute 
\begin{equation} \label{dimberegn}
  \delta^*(Z:\mathcal D) = 2 - \lim_{t\to 0^+}
  (t\Phi^*(D+cT_1+\sqrt{t}Y,(D+cT_1+\sqrt{t}Y)^*:\mathcal D)),  
\end{equation}
where $Y$ is a standard circular element $*$-free from $T_1$, if the limit
exists.
By lemma \ref{distributiontheorem} we know that  
$$
T_1+\sqrt{t}Y\stackrel{*-dist.}{\sim} \sqrt{1+t}T_1+ \sqrt{t}T_2^*,
$$
so using (\ref{scalarfischer}) of remark \ref{fischerremark}
(\ref{dimberegn}) becomes 
\begin{multline*}
    \delta^*(Z:\mathcal D) \\
 =  2 - \lim_{t\to 0^+} (t\Phi^*(D+\sqrt{c^2+t}T_1+\sqrt{t}T_2^*,(D+\sqrt{c^2+t}T_1+\sqrt{t}T_2^*)^*:\mathcal D)) \\
\\  =  2 - \lim_{t\to 0^+} \Phi^*(\tfrac{1}{\sqrt{t}}D+\sqrt{\tfrac{c^2+t}{t}}T_1+T_2^*,(\tfrac{1}{\sqrt{t}}D+\sqrt{\tfrac{c^2+t}{t}}T_1+T_2^*)^*:\mathcal D)
\end{multline*}
if the limit on the right-hand side exists.
Define
\begin{equation*}
  (S_t,S_t^*) = \left(\tfrac{1}{\sqrt{t}}D+\sqrt{\tfrac{c^2+t}{t}}T_1 +
  T_2^*,\left(\tfrac{1}{\sqrt{t}}D+\sqrt{\tfrac{c^2+t}{t}}T_1 + T_2^*\right)^*\right). 
\end{equation*}
Then it easy to see that 
\begin{equation*}
  (\xi_t,\xi_t^*) = \left(\sqrt{\tfrac{t}{c^2+t}}T_1^* +
  T_2,\sqrt{\tfrac{t}{c^2+t}}T_1 + T_2^*\right) 
\end{equation*}
satisfies the conjugate relations for $(S_t,S_t^*)$ with respect to
$\mathcal D$ because $(\xi_t,\xi^*_t,S_t,S_t^*)$ is a (non-centered) $\mathcal
D$-Gaussian system. For example 
\begin{eqnarray*}
  \kappa^\mathcal D_2(\xi_t\otimes_\mathcal D d S_t) & = &  E_\mathcal D\left(\left(\sqrt{\tfrac{t}{c^2+t}}T_1^* + T_2\right)d\left(\sqrt{\tfrac{c^2+t}{t}}T_1 + T_2^*\right)\right) \\
 & = & L^*(d) + L(d) = \left(\int_0^1 d(y)\text{d} y \right)1 = \tau(d)1,
\end{eqnarray*}
for $d\in \mathcal D$ because $T_1$ and $T_2$ are $\mathcal D$-free, and $D\in \mathcal D$ is
$E_\mathcal D$-free from everything \cite[Lemma 3.2.4]{AMS}. The other
$\kappa^\mathcal 
D_2$-identities are checked similarly. Furthermore $\xi_t,\xi_t^*\in
L^2(\{S_t,S_t^*\}\cup \mathcal D,\tau)$ because lemma
\ref{cutouttheorem} says that actually $T_1,T_2\in  W^*(S_t, \mathcal
D)$. 

We thus conclude that
\begin{eqnarray*}
  \Phi^*(S_t,S_t^*:\mathcal D) &  = & ||\xi_t ||_\tau^2 + ||\xi_t^* ||_\tau^2 = 2||\xi_t ||_\tau^2 \\
 &  = & 2\tfrac{t}{c^2+t}\tau(T_1^*T_1)+ 2\tau(T_2T_2^*) \\ 
 & & +2\sqrt{\tfrac{t}{c^2+t}}(\tau(E_\mathcal D(T_1^*T^*_2) + E_\mathcal D(T_2T_1))) \\
& = &  \tfrac{t}{c^2+t} + 1,
\end{eqnarray*}
since $E_\mathcal D$-freeness of $T_1$ and $T_2$ implies that 
$$
E_\mathcal D(T_1^*T_2^*) = E_\mathcal D(T_2T_1)=0.
$$

We have now shown that $(\xi_t,\xi_t^*)$ is a conjugate system for $(S_t,S_t^*)$, so by (\ref{dimberegn})
\begin{equation*}
  \delta^*(Z:\mathcal D) = 2 - \lim_{t\to 0^+} (\tfrac{t}{c^2+t} + 1) = 1. 
\end{equation*}
\qed


\section{$\delta^*$ for a $\DT$-operator with respect to $\mathbb C$}

Let $T_1=\mathcal{UT}(X_1,\lambda)\in \mathcal M$ and
$T_2=\mathcal{UT}(X_2,\lambda)\in \mathcal M$ be two $E_\mathcal D$-free
upper triangular operators and let $D\in \mathcal D$ be chosen such that
$D\in W^*(D+cT_1)$. As in the previous section we define $S_t= 
\tfrac{1}{\sqrt{t}}D+\sqrt{\tfrac{c^2+t}{t}}T_1+T_2^*$, and
$\xi_t=\sqrt{\tfrac{t}{c^2+t}}T_1^*+T_2$. 
$Z=D+cT_1$ is a $\DT$-operator and by lemma \ref{lemma6.2} every $\DT$-operator can be
realized in this way. 

The goal is to compute
$\delta^*(Z) := \delta^*(Z:\mathbb C)$. Inspecting the proof of
$\delta^*(Z:\mathcal D)=1$ from
section \ref{deltaD}  we observe that a
crucial point in showing that $(\xi_t,\xi_t^*)$ is a conjugate system
for $(S_t,S_t^*)$ is the use of the algebra $\mathcal D$  to ``cut out''
the upper and lower triangular parts of $S_t$ and $S_t^*$.  
Since $\mathcal D$ can not be contained in $L^2(\{S_t,S_t^*\},\tau)$
\cite[Prop. p. 123]{Voi6} the
above argument does not apply to compute $\delta^*(Z)$. By remark
\ref{fischerremark} (c) it does however tell us that $(\xi_t,\xi_t^*)$
satisfies the conjugate relations for $(S_t,S_t^*)$ with respect to
$\mathbb C$. Summarizing we have  
\begin{equation}
  \label{summarizedelta}
  \delta^*(Z) \geq 2 - \limsup_{t\to 0^+} \Phi^*(S_t,S_t^*:\mathbb C) \geq 2 - \limsup_{t\to 0^+}\Phi^*(S_t,S_t^*:\mathcal D) = 1
\end{equation}

We would like to compute the Hilbert space projection
$$
P:L^2(\{S_t,S_t^*\}\cup \mathcal D,\tau)\to
L^2(\{S_t,S_t^*\},\tau).
$$ 
This however seems to be very
difficult. Instead our strategy will be to estimate the distance from
$\xi_t$ to a suitable subspace of the orthogonal complement of 
$L^2(\{S_t,S_t\},\tau)$ in
$L^2(\{S_t,S_t^*\}\cup \mathcal D,\tau)$ as $t \to 0^+$. To do this we
first produce an element in the orthogonal complement
$(L^2(\{S_t,S_t^*\},\tau))^\perp$.  
\begin{lemma}  \label{liberationlemma}
\begin{multline} \label{liberationeq} 
  j_t:=[\xi_t,S_t]+[\xi_t^*,S_t^*] \\ = 
\left(\sqrt{\tfrac{c^2+t}{t}}- \sqrt{\tfrac{t}{c²+t}}\right)\left([T_1,T_2]+[T_1^*,T_2^*]\right) \\ 
+ \left[\tfrac{1}{\sqrt{t}}D,\sqrt{\tfrac{t}{c^2+t}}T^*_1+T_2\right]
+\left[\tfrac{1}{\sqrt{t}}D^*,\sqrt{\tfrac{t}{c^2+t}}T_1+T^*_2\right]
\end{multline}
belongs to $(L^2(\{S_t,S_t^*\},\tau))^\perp$ in $L^2(\{S_t,S_t^*\}\cup \mathcal D,\tau)$.
\end{lemma}
\begin{remark}
The above element, $j_t$, is actually just Voiculescu's
liberation gradient of $(\alg(\{S_t,S_t^*\}),\mathcal
D)$ \cite[section 5]{Voi6}. Voiculescu shows that the liberation
gradient $j(A,B)$, of two 
unital algebras $A$ and $B$ in a $W^*$-probability space $(M,\tau)$
measures how $\tau$-free $A$ and $B$ are, in the 
sense that $j(A,B)=0$ if and only if $A$ and $B$ are $\tau$-free. The
liberation gradient $j(A,B)$ satisfies
\begin{multline} \label{liberation}
\tau(j(A,B)a_1b_1\cdots a_nb_n) = \sum_{k=1}^n (\tau(a_1b_1\cdots
b_{k-1}a_k)\tau(b_ka_{k+1}b_{k+1}\cdots a_nb_n) \\
  - \tau(a_1b_1\cdots a_{k-1}b_{k-1})\tau(a_kb_k\cdots a_nb_n)),
\end{multline}
for all $a_1,\ldots,a_n\in A$ and $b_1\ldots,b_n\in B$.
\end{remark}

\begin{proof}
Since $(\xi_t,\xi_t^*)$ is a conjugate system for
$(S_t,S_t^*)$ with respect to $\mathcal D$ it satisfies the conjugate
relations  (\ref{conjugate}). Using (\ref{conjugate}) it is easy to show
by direct computation that $j_t:=
[\xi_t,S_t]+[\xi_t^*,S_t^*]$ satisfies
(\ref{liberation}), that is, 
\begin{multline} \label{liberationD}
\tau(j_t S^{i_1}_t d_1\cdots S_t^{i_n} d_n) = \sum_{k=1}^n 
 (\tau(S_t^{i_1} d_1\cdots
d_{k-1} S_t^{i_{k}})\tau(d_k S_t^{i_{k+1}}d_{k+1}\cdots S_t^{i_n}d_n) \\
  - \tau(S_t^{i_1}d_1\cdots S_t^{i_{k-1}}d_{k-1})\tau(S_t^{i_k}d_k\cdots S_t^{i_n}d_n)),
\end{multline}
for all $d_1,\ldots,d_n\in \mathcal D$ and $i_1,\ldots i_n\in
\{1,*\}$. Restricting (\ref{liberationD}) to the case where $d_1=\cdots
=d_n=1\in \mathcal D$ we infer that 
\begin{equation*}
\tau(j_t S^{i_1}_t \cdots S_t^{i_n}) = 0 
\end{equation*} 
for all $i_1,\ldots,i_n\in \{1,*\}$. The inner product of $j_t$ with
monomials of the form $ S^{i_1}_t \cdots S_t^{i_n}$ is thus 
zero, and since these monomials span a dense linear subspace of 
$L^2(\{S_t,S_t^*\},\tau)$ we conclude that $j_t\perp L^2(
\{S_t,S_t^*\},\tau)$. The last equality of
(\ref{liberationeq}) is easily checked by direct computation.
\end{proof}

We now use $j_t\in W^*(S_t)^\perp$ from (\ref{liberationeq}) to give a
new lower bound of $\delta^*(Z)$.

\begin{lemma} \label{commutatorlemma}
  \begin{multline*}
\limsup_{t\to\ 0^+}\Phi^*(S_t,S_t^*:\mathbb C) \\ \leq
2\dist_2(T_2,W^*(T_1)([D+cT_1,T_2]+[(D+cT_1)^*,T_2^*])W^*(T_1))^2
\end{multline*}
\end{lemma}
\begin{proof}
Define $\tilde{S}_t = \sqrt{t}S_t =
D+\sqrt{c^2+t}T_1 + \sqrt{t}T_2^*$. We know that 
\begin{equation} \label{tildealpha}
\Phi^*(S_t,S_t^*:\mathbb C) =
2||E_{W^*(\tilde{S}_t)}(\xi_t) ||_\tau^2 
\end{equation}
Let $A\in W^*(D+cT_1)([D+cT_1,T_2]+[(D+cT_1)^*,T_2^*])W^*(D+cT_1)$ and let
$\epsilon>0$. Choose polynomials $p$ and $q$ in $D+cT_1$ and $(D+cT_1)^*$ such
that 
\begin{multline*}
||A- p(D+cT_1,(D+cT_1)^*)([D+cT_1,T_2] 
\\ +[(D+cT_1)^*,T_2^*])q(D+cT_1,(D+cT_1)^*)||_\tau <
\epsilon.
\end{multline*}
Define
\begin{multline*}
B(t) =
p(\tilde{S}_t,\tilde{S}_t^*)\sqrt{t}j_tq(\tilde{S}_t,\tilde{S}_t^*)
 \\
= p(\tilde{S}_t,\tilde{S}_t^*)
 \left(\left(\sqrt{c^2+t}- \tfrac{t}{\sqrt{c²+t}}\right)
([T_1,T_2]+[T_1^*,T_2^*])\right. \\ 
+ \left.[D,\sqrt{\tfrac{t}{c^2+t}}T^*_1+T_2]
+[D^*,\sqrt{\tfrac{t}{c^2+t}}T_1+T^*_2 ]\right)
q(\tilde{S}_t,\tilde{S}_t^*)
\end{multline*}
for $t\geq0$, and observe that
\begin{multline*}
  B(0) = p(D+cT_1,(D+cT_1)^*)([D+cT_1,T_2] \\+[(D+cT_1)^*,T_2^*])
q(D+cT_1,(D+cT_1)^*).
\end{multline*}
Then $||A-B(0)||_\tau <\epsilon$ and
$||B(t)-B(0)|| \to 0$ for $t\to 0^+$.  Since by lemma \ref{liberationlemma}
$j_t\sqrt{t}\in W^*(\tilde{S}_t^*)^\perp$ also $B(t)\in
W^*(\tilde{S}_t^*)^\perp$ so
\begin{equation*}
||E_{W^*(\tilde{S}_t^*)}(A)||_\tau \leq
  ||E_{W^*(\tilde{S}_t)}(B(t))||_\tau + ||A-B(t)||_\tau = 0 + ||A-B(t)||_\tau,
\end{equation*}
and thus
\begin{equation*}
\limsup_{t\to 0^+}||E_{W^*(\tilde{S}_t)}(A)||_\tau \leq
||A-B(t)||_\tau <\epsilon.
\end{equation*}
Since $\epsilon$ was chosen arbitrarily we conclude that
\begin{equation} \label{distarg}
\lim_{t\to 0^+}||E_{W^*(\tilde{S}_t)}(A)||_\tau = 0,
\end{equation}
for all $A\in W^*(D+cT_1)([D+cT_1,T_2]+[(D+cT_1)^*,T^*_2])W^*(D+cT_1)$. 
Now define
\begin{equation*}
\alpha = \dist_2(T_2,W^*(T_1)([D+cT_1,T_2]+[(D+cT_1)^*,T_2^*])W^*(T_1)).
\end{equation*}
Since $D\in \mathcal D$ is chosen such that $D\in W^*(D+cT_1)$ we have 
$W^*(T_1)=W^*(D+cT_1)$ so from (\ref{distarg}) we conclude that
\begin{equation*}
\limsup_{t\to 0^+} ||E_{W^*(\tilde{S}_t)}(T_2)||_\tau \leq \alpha,
\end{equation*}
but since $||E_{W^*(\tilde{S}_t)}||_\tau\leq 1$ and since
$||\xi_t-T_2||_\tau \to 0$ for $t\to 0^+$ we conclude that 
\begin{equation}
\limsup_{t\to
0^+}||E_{W^*(\tilde{S}_t)}(\xi_t)||_\tau \leq \alpha.
\end{equation}
Combining with (\ref{tildealpha}) we have
\begin{equation*}
\limsup_{t\to 0^+}\Phi^*(S_t,S_t^*:\mathbb C) \leq 2\alpha^2.
\end{equation*}
\end{proof}

Lemma \ref{commutator} can now be used to get rid of the
$[(D+cT_1)^*,T_2^*]$-term in lemma  \ref{commutatorlemma}.

\begin{proposition} \label{commutatortheorem}
  \begin{equation}
\limsup_{t\to\ 0^+}\Phi^*(S_t,S_t^*:\mathbb C) \leq
2\dist_2(T_2,W^*(T_1)([D+cT_1,T_2])W^*(T_1))^2.
\end{equation}
\end{proposition}

\begin{proof}
It is immediate that
$$
E=\overline{W^*(T_1)([D+cT_1,T_2]+[(D+cT_1)^*,T_2^*])W^*(T_1)}^{||\cdot
||_\tau}
$$
 is invariant under multiplication from the left and right with
elements from $W^*(T_1)$.
By lemma \ref{commutator} we know that $[D+cT_1,T_2]\in E$ so
\begin{equation*}
W^*(T_1)[D+cT_1,T_2]W^*(T_1)\subseteq E.
\end{equation*} 
We thus have 
\begin{multline*}
\dist_2(T_2,W^*(T_1)[D+cT_1,T_2]W^*(T_1))  \geq  \dist_2(T_2,E) \\
 =  \dist_2(T_2,W^*(T_1)([D+cT_1,T_2]+[(D+cT_1)^*,T_2^*])W^*(T_1)),
\end{multline*}
so the proposition follows from lemma \ref{commutatorlemma}. 
\end{proof}

Let $W^*(T_1)^0$ be the opposite algebra, that is, $W^*(T_1)^0$ is
just $W^*(T_1)$ with multiplication reversed. 

\begin{lemma} \label{statelemma}
Let $\Hat{\otimes}$ denote the von Neumann algebra tensor product, and equip
$W^*(T_1)$ with the usual faithful normal tracial state, $\tau$, and usual
conditional expectation, $E_\mathcal D$, given by restriction to $W^*(T_1)$.
There exists a positive functional, $\phi: W^*(T_1) \Hat{\otimes}
W^*(T_1)^0\to \mathbb 
C$ such that  
\begin{enumerate}
\item[(i)] $0 \leq \phi \leq \tau\otimes \tau^0$ 
\item[(ii)] $\tau(T_2^*a T_2b) = \phi (a \otimes b^0)$ for $a,b\in W^*(T_1)$.
\end{enumerate}
where $b\mapsto b^0:W^*(T_1)\to W^*(T_1)^0$ is the anti-multiplicative
isomorphism. 
\end{lemma}
\begin{proof}
For all $a,b\in W^*(T_1)$ we have
\begin{eqnarray*}
\tau(T_2^*aT_2b) & = & \int_0^1\left(\int_0^x E_\mathcal
D(a)(t)\text{d}t\right)E_\mathcal D(b)(x)\text{d} x \\
& = & \int_0^1\int_0^1 E_\mathcal D(a)(t)E_\mathcal D(b)(x)h(t,x)\text{d} t
\text{d} x,
\end{eqnarray*}
where
\begin{equation*}
h(t,x) = 
\begin{cases}
1, & t\leq x \\
0, & t > x
\end{cases}
\end{equation*}
for $0\leq x \leq 1$.
$h(t,x)$ corresponds to an element $H\in D\Hat{\otimes} D^0$ such that
$0\leq H\leq 1$, so 
\begin{equation*}
\tau(T_2^*aT_2b)= (\tau\otimes \tau^0)(H(E_\mathcal D(a)\otimes E_{\mathcal D^0}(b))).
\end{equation*}
Since $E_\mathcal D \otimes E_{\mathcal D^0}$ is a positive normal
operator on $W^*(T_1)\hat{\otimes}W^*(T_1)^0$ we have
\begin{equation*}
\phi:z\mapsto (\tau\otimes \tau^0)(H(E_\mathcal D\otimes E_{\mathcal D^0}(z)))
\end{equation*}
is a positive normal functional on $W^*(T_1)\Hat{\otimes}W^*(T_1)^0$. 
For $z\geq 0$ we observe that
\begin{equation*}
\phi(z) \leq (\tau\otimes \tau^0)((E_\mathcal D\otimes E_{\mathcal
D^0})(z)) = (\tau\otimes \tau^0)(z),
\end{equation*}
so
\begin{equation*}
0\leq \phi\leq \tau\otimes \tau^0.
\end{equation*}
\end{proof}

We want to estimate the distance $\alpha=
\dist_2(T_2,W^*(T_1)[D+cT_1,T_2]W^*(T_1))$ from proposition
\ref{commutatortheorem}. Let $\phi:W^*(T_1)\Hat{\otimes}W^*(T_1)^0\to
\mathbb C$ be the state from lemma \ref{statelemma}.
We observe that for $a_1,\ldots,a_n,b_1,\ldots,b_n\in W^*(T_1)$ lemma
\ref{statelemma} implies that
\begin{eqnarray*}
||T_2-\sum_{i=1}^n a_iT_2 b_i||_\tau^2 & = & ||T_2||_\tau^2 -2\Re
  \tau\left(\sum_{i=1}^n T_2^*a_iT_2b_i \right) \\
 &  &  +
  \tau\left(\sum_{i,j=1}^n b_i^*T_2^*a_i^*a_jT_2b_j \right) \\
 & = & \phi(1\otimes 1) - 2\Re \phi\left( \sum_{i=1}^n a_i\otimes b_i^0
  \right) \\
 & & + \phi\left( \sum_{i,j=1}^n a_i^*a_j\otimes (b_jb_i^*)^0\right) \\
 & = & \phi\left(\left(1\otimes 1 - \sum_{i=1}^n a_i\otimes b_i^0
  \right)^* \left(1\otimes -\sum_{j=1}^n a_j\otimes b_j^0
  \right)\right) \\
 & \leq & \left|\left |1\otimes 1 - \sum_{i=1}^n a_i\otimes
  b_i^0\right|\right|^2_{\tau\otimes \tau^0}.
\end{eqnarray*}
We thus have
\begin{eqnarray*}
\alpha & = &
\dist_2(T_2,W^*(T_1)[D+cT_1,T_2]W^*(T_1)) \\
 & = & \inf_{\substack{n\in \mathbb N \\ a_i,b_i\in W^*(T_1)}} \left| \left|
 T_2 - \sum_{i=1}^n a_i((D+cT_1)T_2-T_2(D+cT_1))b_i \right|\right|_\tau \\
& = & \inf_{\substack{n\in \mathbb N \\ a_i,b_i\in W^*(T_1)}}
\left|\left| T_2 - \sum_{i=1}^n (a_i(D+cT_1))T_2b_i - a_iT_2((D+cT_1)b_i)
\right|\right|_\tau \\
& \leq & \inf_{\substack{n\in \mathbb N \\ a_i,b_i\in W^*(T_1)}}
\left|\left| 1\otimes 1 - \sum_{i=1}^n (a_i(D+cT_1)\otimes b_i^0 - a_i\otimes
b_i^0(D+cT_1)^0) \right|\right|_{\tau\otimes \tau^0} \\
& = & \inf_{\substack{n\in \mathbb N \\ a_i,b_i\in W^*(T_1)}}
\left|\left|  1\otimes 1 - \left(\sum_{i=1}^n (a_i\otimes
b_i^0)\right)((D+cT_1)\otimes 1 - 1\otimes (D+cT_1)^0) \right| \right|_{\tau\otimes
\tau^0} .
\end{eqnarray*}
If we can show that $\ker(R)=\{0\}$ for $R:=(D+cT_1)\otimes 1 - 1\otimes
(D+cT_1)^0$ then defining $Q_n=f_n(R^*R)R^*$ where 
\begin{equation*}
f_n(t) = 
\begin{cases}
0 & \textrm{for } 0\leq t \leq \tfrac{1}{n} \\
\tfrac{1}{t} & \textrm{for } t \geq \tfrac{1}{n},
\end{cases}
\end{equation*}
we have $Q_nR = 1_{[\tfrac{1}{n},\infty)}(R^*R)$ so since
$\ker(R^*R)=\ker(R)=\{0\}$ we have $\norm{Q_nR-1\otimes 1}_{\tau\otimes
\tau^0}\to 0$. Since each $Q_n\in W^*(T_1)\Hat{\otimes}W^*(T_1)^0$
can be approximated in the strong operator topology by operators of the
form $\sum_{i=1}^m a_i\otimes b_i^0$ where $a_i\in W^*(T_1)$ and $b_i^0\in
W^*(T_1)^0$ it follows that $\alpha=0$. 

 Combining this with
theorem \ref{commutatortheorem} we have
\begin{equation}
 \ker((D+cT_1)\otimes 1 - 1\otimes (D+cT_1)^0)=\{0\}\Rightarrow
 \limsup_{t\to 0^+} \Phi^*(S_t,S_t^*:\mathbb C)=0. 
\end{equation}
By (\ref{summarizedelta}) this will imply that $\delta^*(T)\geq 2$. So
now the only remaining problem is to show that $ \ker((D+cT_1)\otimes 1 -
1\otimes (D+cT_1)^0)=\{0\}$, and this will follow from the following
``eigenspace''-results.

\begin{lemma} \label{Mondelemma}
Let $A,B\in \mathbb B(H)$ be bounded operators on a Hilbert space, $H$,
such that $\ker(B)=\{0\}$. Define $E_\lambda = \{x\in H| Ax=\lambda Bx\}$. If
$\lambda_1,\ldots,\lambda_n\in \mathbb C$ are mutually different, then 
the corresponding subspaces $E_{\lambda_1},\ldots,E_{\lambda_n}$ are all
linearly independent.  
\end{lemma}

\begin{proof}
Let $\lambda_1,\ldots,\lambda_n\in \mathbb C$ be mutually different
complex numbers. We must show that $\sum_{i=1}^n x_i=0$ implies
$x_1=\cdots =x_n=0$ when $x_i\in E_{\lambda_i}$ for all $i\in
\{1,\ldots,n\}$. But $\sum_{i=1}^n x_i=0$ implies that
\begin{equation*}
\sum_{i=1}^n A^kx_i=0,
\end{equation*} 
for all $k\in \{0,\ldots,n-1\}$. Since $x_i\in E_{\lambda_i}$ we have
\begin{equation*}
0=\sum_{i=1}^n \lambda_i^kB^kx_i.
\end{equation*}
Multiplying by $B$ an appropriate number of times we obtain
\begin{equation*}
  \begin{bmatrix}
  1 & 1 & \cdots & 1\\
 \lambda_1 & \lambda_2 & \cdots & \lambda_n \\
 \vdots & \vdots & \ddots & \vdots \\
\lambda_1^{n-1} & \lambda_2^{n-1} & \cdots & \lambda_n^{n-1}
  \end{bmatrix}
  \begin{bmatrix}
  B^{n-1}x_1 \\
  B^{n-1}x_2 \\
  \vdots \\
  B^{n-1}x_n 
  \end{bmatrix} =0.
\end{equation*}
Since the determinant on the left hand-side is a van der
Monde determinant, which is exactly invertible for $\lambda_1,\ldots
,\lambda_n$ all mutually different, we infer that $B^{n-1}x_1=\cdots
B^{n-1}x_n=0$. Since $\ker(B)=\{0\}$ also $\ker(B^{n-1})=\{0\}$ so
$x_1=\cdots=x_n=0$. Thus the ``eigenspaces'' 
$E_{\lambda_1},\ldots,E_{\lambda_n}$ are linearly independent.

\end{proof}
\begin{lemma} \label{eigenvalue}
Let $\mathcal N\subseteq \mathbb B(H)$ be a finite $W^*$-algebra represented on a Hilbert space
$H$. Let $A,B\in \mathcal N$, and
define $E_\lambda = \{x\in H| Ax=\lambda Bx\}$. If $\ker B=\{0\}$ then  
\begin{equation}
E_{\lambda_n} \cap \overline{E_{\lambda_1}+\cdots + E_{\lambda_{n-1}}} = \{0\},
\end{equation}
when $\lambda_1,\ldots, \lambda_n\in \mathbb C$ are mutually different.
\end{lemma}

\begin{proof}
By lemma \ref{Mondelemma} we know that
$E_{\lambda_1},\ldots,E_{\lambda_n}$ are linearly independent for
$\lambda_1,\ldots,\lambda_n$ mutually different, so we have
\begin{equation}
E_{\lambda_n}\cap\left(E_{\lambda_1}+\cdots + E_{\lambda_{n-1}}\right) =\{0\}.
\end{equation}
Obviously $E_{\lambda_1}+\cdots + E_{\lambda_{n-1}}$ and
$E_{\lambda_n}$  are subspaces of $H$
affiliated with $\mathcal N$, so theorem \ref{affsubspaces} implies that
\begin{eqnarray*}
E_{\lambda_n}\cap \overline{E_{\lambda_1}+\cdots+E_{\lambda_{n-1}}} &
\subseteq & \overline{E_{\lambda_n}}\cap
\overline{E_{\lambda_1}+\cdots+E_{\lambda_{n-1}}} \\
& = & \overline{E_{\lambda_n}\cap
\left(E_{\lambda_1}+\cdots+E_{\lambda_{n-1}}\right)}=\{0\}.
\end{eqnarray*}
\end{proof}

\begin{proposition} \label{kaplansky}
Let $A,B\in \mathbb B(H)$ be bounded operators on a Hilbert space, $H$,
such that $\ker B=\{0\}$  and
define $E_\lambda = \{x\in H| Ax=\lambda Bx\}$. Assume that $A,B\in
\mathcal N$
where $\mathcal N$ is a $II_1$-factor with a faithful tracial state, $\tau$. Let
$p_\lambda$ be 
the projection onto $E_\lambda$. If $\lambda_1,\ldots,\lambda_n$ are
mutually different, then 
\begin{equation}
\tau\left(\bigvee_{j=1}^n p_{\lambda_j} \right) = \sum_{j=1}^n
\tau(p_{\lambda_j}). 
\end{equation}
\end{proposition}

\begin{proof}
If $p,q$ are two projections in $\mathcal N$ then it follows
from Kaplansky's formula that 
\begin{equation} \label{kaplformula}
\tau(p\vee q) = \tau(p) + \tau(q) -\tau(p\wedge q).
\end{equation}
Define $q_k = \bigvee_{i=1}^k p_j$ for $k=1,\ldots,n$, and put
$q_0=0$. Then $q_k=p_k \vee q_{k-1}$ and by  lemma \ref{eigenvalue}
$p_k\wedge q_{k-1}=0$ for $k\geq 1$
since $E_{\lambda_k} \cap
\overline{E_{\lambda_1}+\cdots 
+E_{\lambda_{k-1}}}=\{0\}$. We thus conclude from (\ref{kaplformula}) that
\begin{equation*}
\tau(q_k)= \tau(p_k)+ \tau(q_{k-1}),
\end{equation*}
and thus also
\begin{equation*}
\tau(q_n) = \sum_{j=1}^n\tau(p_j).
\end{equation*}
\end{proof}
 
\begin{corollary} \label{countablecor}
Let $A,B\in \mathbb B(H)$ be bounded operators on a Hilbert space, $H$,
such that $\ker B=\{0\}$  and
define $E_\lambda = \{x\in H| Ax=\lambda Bx\}$. Assume that $A,B\in
\mathcal N$
where $\mathcal N$ is a $II_1$-factor with a faithful tracial state, $\tau$. Let
$p_\lambda$ be 
the projection onto $E_\lambda$. Then $E_\lambda = \{0\}$ except for
countably many $\lambda\in \mathbb C$.
\end{corollary}

\begin{proof}
For any finite subset $F\subset \mathbb C$ proposition
\ref{kaplansky} implies that
\begin{equation*}
\sum_{\lambda\in F} \tau(p_\lambda) = \tau\left(\bigvee_{\lambda\in
F}p_\lambda \right) \leq 1,
\end{equation*}
so
\begin{equation*}
\sum_{\lambda\in \mathbb C} \tau(p_\lambda)\leq 1,
\end{equation*}
and thus $\tau(p_\lambda)=0$ for all but countably many $\lambda \in
\mathbb C$.
\end{proof}

\begin{corollary} \label{corker}
  \begin{equation*}
  \ker((D+cT_1)\otimes 1 - 1\otimes (D+cT_1)^0)=\{0\}.
  \end{equation*}
\end{corollary} 
\begin{proof}
Define $A_1= T_1\otimes 1$ and $B_1= 1\otimes T_1^0$. It follows from
\cite[Th. 8.9]{DT} that $\ker(B_1)=\ker(1\otimes T_1^0)=\{0\}$ so corollary
\ref{countablecor} applies. $\ker(T_1\otimes
1-1\otimes T_1^0)\neq \{0\}$ would imply that $\ker(T_1\otimes
1-1\otimes \lambda T_1^0)\neq \{0\}$ for any $\lambda \in \mathbb T$ in
the unit 
circle of $\mathbb C$, since $T_1$ and $\lambda T_1$ have the same
$*$-distribution \cite[Prop. 2.12]{DT}. But then $E_\lambda\neq \{0\}$ for all $\lambda
\in \mathbb T$ which is not countable. This contradicts corollary
\ref{countablecor} so we must have $\ker(T_1\otimes 1 - 1\otimes
T_1^0)=\{0\}$.
Next define $A=D\otimes 1 - 1\otimes D^0$ and $B=T_1\otimes 1 - 1\otimes
T_1^0$. Then we have just seen that $\ker(B)=\{0\}$, so since $(A,B)$
and $(A,\lambda B)$ are equally $*$-distributed for any $\lambda\in
\mathbb T$ we can use corollary \ref{countablecor} once more to conclude
that $\ker((D+cT)\otimes 1 - 1\otimes (D+cT)^0)=\{0\}$.
\end{proof}

Corollary \ref{corker} finishes the proof of $\delta^*(T)\geq 2$. It
follows from proposition 
\ref{deltaovre} used on the real and imaginary part of $T$ that
$\delta^*(T)\leq 2$ so we conclude that $\delta^*(T)=2$. 

Proceeding as in the proof of corollary \ref{corker} we have
following corollary on the point spectrum of $\DT$-operators.
\begin{corollary}
Every $\DT$-operator has empty point spectrum.
\end{corollary}
\begin{proof}
Let $D+cT$ be a $\DT$-operator for $D\in \mathcal D$ and $T\in \mathcal{UT}(X,c)$ for some compactly supported complex probability measure, $\mu$,
 and some $c>0$. Let $\gamma\in \mathbb C$ be fixed. 
Then $(-\gamma 1 +D)+cT$ is again a $\DT$-operator and the $*$-distribution of $(-\gamma 1+ D)+cT$ is completely determined by the distribution of 
$(-\gamma 1+D)$ and $c$.
Since $\lambda T$ and $T$ are equally $*$-distributed for all $\lambda\in \mathbb T$ we infer that $(-\gamma 1+D)+cT$ and $(-\gamma 1 +D)+\lambda cT$ are equally $*$-distributed for all $\lambda\in \mathbb T$. By an argument similar to 
the one given in corollary \ref{corker} we infer that 
\begin{equation*}
  \ker(\gamma 1 - (D+cT))=\ker((-\gamma 1+D)+cT)=\{0\}
\end{equation*}
for all $\gamma\in \mathbb C$, so $\sigma_p(D+cT)=\emptyset$.

\end{proof}

\bibliographystyle{amsplain}

\begin{thebibliography}{25}
 

\bibitem[AZ]{AZ} P. Ainsworth, \emph{Ubegrænsede operatorer affilieret
med en endelig von Neumann algebra}, Master Thesis, University of Odense, 1985.

\bibitem[DH1]{DT} K. Dykema, U. Haagerup, \emph{DT-operators and
decomposability of Voiculescu's circular operator}, preprint. To appear
in Amer. J. Math.

\bibitem[DH2]{invsub} \bysame, \emph{Invariant subspaces of the
quasinilpotent DT-operator}, preprint. To appear in JFA.

\bibitem[HP]{HiaiPetz} F. Hiai, D. Petz, \emph{The Semicircle Law, Free
random Variables and Entropy}, Mathematical Surveys and Monographs,
Vol. 77, Am Math. Soc., 2000.

\bibitem[NSS1]{NSS1} A. Nica, D. Shlyakhtenko, R. Speicher, \emph{Some
minimization problems for the free analogue of the free Fischer
information}, Advances in Mathematics \textbf{141} (1999), no. 2, 282-321.

\bibitem[NSS2]{NSS2} \bysame, \emph{Operator-valued
distributions. I. Characterizations of freeness},
Int. Math. Res. Not. \textbf{2002}, no. 29, 1509-1538.

\bibitem[NSS3]{NSS3} \bysame, \emph{A Characterization of Freeness by a
Factorization Property of R-transforms}, MSRI Preprint 2001-001.

\bibitem[KR]{KR} R.V. Kadison, J.R. Ringrose, \emph{Fundamentals of the Theory of Operator Algebras}, Vol. II, Academic Press, New York, 1983, 1986. 

\bibitem[Sk]{Skau} C. Skau, \emph{Finite subalgebras of a von Neumann
algebra}, Journal of Functional analysis \textbf{25} (1977), 211-235.

\bibitem[\'Sn]{SniadyDT} P. \'Sniady, \emph{Inequality for Voiculescu's
entropy in terms of Brown measure},
Internat. Math. Res. Notices \textbf{2003}, 51-64. 


\bibitem[\'SnSp]{SnSp} P. \'Sniady, R. Speicher, \emph{Continuous family
of invariant subspaces for R-diagonal operators},
Invent. Math. \textbf{146} (2001), no. 2, 329-363.

\bibitem[Sp1]{AMS} R. Speicher, \emph{Combinatorial theory of the free product with amalgamation and operator-valued free probability theory}, Mem. Amer. Math. Soc. \textbf{132} (1998), no 627, x+88 pp. 

\bibitem[Sp2]{Combinatorics} R.Speicher, \emph{Combinatorics of free probability
theory}, ``Free probability and operator spaces'', IHP, Paris, 1999.

\bibitem[Ta]{Ta} M. Takesaki \emph{Theory of operator algebras I},
Encyclopaedia of Mathematical Sciences Vol. 124, Operator Algebras and
Non-commutative geometry Vol. 5, Springer-Verlag, Berlin Heidelberg New
York, 2000.

\bibitem[Voi1]{Voi1} D. Voiculescu, \emph{The analogues of entropy and
of Fischer's information measure in free probability theory, I}
Communications Math. Physics \textbf{155} (1993), 71-92.

\bibitem[Voi2]{Voi2} \bysame,  \emph{The analogues of entropy and
of Fischer's information measure in free probability theory, II},
Invent. Math. \textbf{118} (1994) 411-440.

\bibitem[Voi5]{Voi5} \bysame, \emph{The analogues of entropy and
of Fischer's information measure in free probability theory, V:
Non-commutative Hilbert Transforms}, Invent. Math. \textbf{132}
(1998), 189-227.

\bibitem[Voi6]{Voi6} \bysame, \emph{The analogues of entropy and
of Fischer's information measure in free probability theory, VI:
Liberation and Mutual Free Information}, Advances in Mathematics
\textbf{146} (1999), 101-166.

\end{thebibliography}
 \providecommand{\bysame}{\leavevmode\hbox to3em{\hrulefill}\thinspace}

\appendix
\section{Affiliated subspaces}

This appendix is part of an unpublished master thesis \cite{AZ} done in 1985 by
P. Ainsworth under supervision of U. Haagerup. The translated title of
the thesis is: \emph{Unbounded operators 
affiliated with a finite von Neumann algebra}. 

\begin{definition}
Let $\mathcal A$ be a finite von Neumann algebra represented on a 
Hilbertspace, $H$. Let $\mathcal E$ be a subspace of $H$. We say that
$\mathcal E$ is \emph{affiliated} to $\mathcal A$ if for all $A'$ in the
commutant, $\mathcal A'$, of $\mathcal A$, and for all $\xi\in \mathcal E$ we have $A'\xi\in
\mathcal E$. 

We say that an operator, $T$, is affiliated with the von Neumann
algebra, $\mathcal A$ if $T$ for every $A'\in \mathcal A'$ we have
$A'\mathscr D(T) \subseteq \mathscr D(T)$ and for all $\xi \in
\mathscr D(T)$ we have $A'T\xi=TA'\xi$, i.e. if $TA'\subset A'T$.
\end{definition}
It is an easy exercise to see that closures and intersections of
affiliated subspaces are again affiliated subspaces.

It is the purpose of this appendix is to give a proof of the following
theorem.   
\begin{theorem} \label{affsubspaces}
Let $\mathcal A$ be a finite von Neumann algebra represented on a Hilbert space
$H$. Let  $\mathcal E$ and
$\mathcal F$ be subspaces of $H$ affiliated with $\mathcal A$. Then 
\begin{equation*}
\overline{\mathcal E}\cap \overline{\mathcal F} = \overline{\mathcal
E\cap \mathcal F}.
\end{equation*}
\end{theorem}

The proof relies on the $T$-theorem.
\begin{theorem}[T-theorem]\cite[Cor. 2]{Skau}
Let $\mathcal A$ be  a finite von Neumann algebra on a Hilbert space
$H$ and let $\xi\in H$. If $\eta\in \overline{\mathcal A\xi}$ there
exists a closed densely defined operator, $T$, affiliated with
$\mathcal A$ such that $\xi\in \mathscr D(T)$ and $T\xi=\eta$.
\end{theorem}  

Actually we need the following version of the T-theorem for the
commutant of a finite von Neumann algebra. 

\begin{corollary} \label{TA'}
Let $\mathcal A$ be  a finite von Neumann algebra on a Hilbert space
$H$. Then the $T$-theorem is valid for the commutant of $\mathcal A$,
that is for all $\eta \in \overline{\mathcal A'\xi}$ there 
exists a closed densely defined operator, $T$, affiliated with
$\mathcal A'$ such that $\xi\in \mathscr D(T)$ and $T\xi=\eta$.
\end{corollary}

\begin{proof}
  Let $\xi,\eta\in H$ such that $\eta \in \overline{\mathcal A'\xi}$ and let $p_\xi$ and $p_\eta$ be the projections onto $\overline{\mathcal A\xi}$ and
$\overline{\mathcal A\eta}$ respectively. By \cite[prop. 9.1.2.]{KR} $p_\xi$ and $p_\eta$ are finite projections in $\mathcal A'$. By \cite[6.3.8.]{KR} the finiteness of $p_\xi$ and $p_\eta$ implies that $p:= p_\xi \vee p_\eta$ is a 
finite projection in $\mathcal A'$. Now $\xi,\eta \in p(H)$, and since the commutant of $\mathcal A\restriction_{p(H)}$ is the finite von Neumann algebra $p\mathcal A'p$ the T-theorem (for the finite $W^*$-algebra $p\mathcal A'p$) says that there exists a closed, densely defined (in $p(H)$) operator, $T$, affiliated with $p\mathcal A'p$ such $T\xi=\eta$. Now extend $T$ to a closed operator which is densely defined in $H$ be defining $T=0$ on the orthogonal complement of $p(H)$. Then $T$ is a closed, densely defined operator affiliated with $\mathcal A'$ such that $\xi\in \mathscr D(T)$ and $T\xi =\eta$.    
\end{proof}

\begin{lemma} \label{projectionlemma}
Let  $\mathcal A$ be  a finite von Neumann algebra on a Hilbert space
$H$. If $(p_i)_{i\in I}$ is an increasing net of projections in
$\mathcal A$ and if $p$ and $q$ are projections in $\mathcal A$ such that $p_i
\nearrow p$ and $q\leq p$ then $p_i\wedge q\nearrow p\wedge q=q$ 
\end{lemma}

\begin{proof}
  By Kaplansky's formula we have
  \begin{equation*}
    q-q\wedge p_i \sim q\vee p_i- p_i.
  \end{equation*}
Since $p_i\leq p$ and $q\leq p$ we have $p_i\vee q\leq p$, so
\begin{equation*}
   q-q\wedge p_i \sim q\vee p_i- p_i \leq p-p_i.
\end{equation*}
$(q\wedge p_i)_{i\in I}$ is an increasing net of projections bounded above by $p$, so $(q\wedge p_i)_{i\in I}$ converges to some projection $r\leq p$.  
For every normal tracial state, $\tau$, we now have
\begin{equation*}
  \tau(q-r)=\lim_i \tau(q-q\wedge p_i) \leq \lim_i \tau(p-p_i)= 0.
\end{equation*}
Since a finite von Neumann algebra has a faithful family of normal tracial states \cite[Th. 2.4.]{Ta} we infer that $q=r$, so $q\wedge p_i\nearrow q$.
\end{proof}

As a final step towards proving theorem \ref{affsubspaces} we need the
following key result. 
\begin{theorem} \label{hjaelpetheorem}
Let $\mathcal A$ be a finite von Neumann algebra represented on a Hilbert space
$H$. Let $\mathcal E$ be a subspace of $H$ affiliated with $\mathcal
A$. Then for all $\xi\in \overline{\mathcal E}$ and for all
$\delta>0$ there exists $A'\in \mathcal A'$ such that $A'\xi\in \mathcal
E$ and $\norm{A'\xi-\xi}<\delta$.
\end{theorem}

\begin{proof}
The proof has three parts. First assume that $\mathcal E= \mathcal
A'\eta$ for some $\eta\in H$. Let $\xi \in \overline{\mathcal E}$. By corollary \ref{TA'} the T-theorem is true for $\mathcal A'$ so there exists a closed 
densely defined operator, $T$, affiliated with $\mathcal A'$ such
that $\eta\in \mathscr D(T)$ and $\xi = T\eta$. Let $T=|T^*|U$ be the
right polar decomposition of 
$T$ and define $p_n=1_{[0,n]}(|T^*|)$. Then $p_n\xi\to \xi$ in norm as
$n\to \infty$. Now
\begin{equation*}
p_n\xi = 1_{[0,n]}(|T^*|)|T^*|U\eta = \left(\int_0^n \lambda
de_\lambda\right)U\eta,
\end{equation*}  
where $e_\lambda$ is the spectral measure of $|T^*|$. Since
$\int_0^n\lambda de_\lambda\in \mathcal A'$ and $U\in
\mathcal A'$ we 
have $p_n\xi \in \mathcal A'\eta=\mathcal E$. Now choose $n$ large enough
to ensure that $\norm{p_n\xi-\xi}<\delta$. This proves the
theorem in the first case.

Secondly assume that $\mathcal E= \sum_{k=1}^n \mathcal A'\eta_k$ for some
$\eta_1,\ldots,\eta_n\in H$. Define
$\tilde{\eta}=(\eta_1,\ldots,\eta_n)\in H^{\oplus n}$. and observe
that for $A_1',\ldots,A_n'\in \mathcal A'$ we have
\begin{equation*}
  \begin{pmatrix}
  A_1' & \cdots & A_n' \\
  0    & \cdots & 0 \\
  \vdots & \ddots & \vdots \\
  0 & \cdots & 0
  \end{pmatrix} 
  \begin{pmatrix}
  \eta_1 \\  \vdots \\ \eta_n 
  \end{pmatrix} = 
  \begin{pmatrix}
  \sum_{k=1}^n A_k'\eta_k \\ 0 \\ \vdots \\ 0
  \end{pmatrix}
\end{equation*}
so $\tilde{\xi}:=(\xi, 0 , \ldots, 0) \in \overline{M_n(\mathcal
A')\tilde{\eta}}$. Identify $\mathcal A$ with its unital embedding as
diagonal 
operators in $M_n(\mathcal A)$. With this identification $\mathcal A$ is
the commutant of $M_n(\mathcal A')$, so the 
idea is now to use the first part of the proof on the subspace
$M_n(\mathcal A')\tilde{\eta}$ of $H^{\oplus n }$ and the von Neumann
algebra $\mathcal A$ represented on $H^{\oplus n}$ with commutant
$M_n(\mathcal A')$. The first part of the proof thus gives us an $A'\in
M_n(\mathcal A')$ 
such that $\norm{A'\tilde{\xi}-\tilde{\xi}}<\delta$, and
$A'\tilde{\xi}\in M_n(\mathcal A')\tilde{\eta}$. If $A'=(A'_{ij})_{i,i=1}^n$ then of
course 
\begin{equation*}
A'\tilde{\xi} = 
\begin{pmatrix}
A_{11}' & \cdots & A_{1n}' \\
\vdots & \ddots & \vdots \\
A_{n1}' & \cdots & A_{nn}' 
\end{pmatrix}
\begin{pmatrix}
\xi \\ 0 \\ \vdots \\ 0
\end{pmatrix}=
\begin{pmatrix}
A_{11}'\xi \\ \vdots \\ A_{n1}'\xi 
\end{pmatrix},
\end{equation*}
so
\begin{equation*}
\delta^2 > \norm{A'\tilde{\xi} - \tilde{\xi}}^2 > \norm{A_{11}'\xi - \xi}^2
\end{equation*}
and furthermore since $A'\tilde{\xi}\in M_n(\mathcal A')\tilde{\eta}$
there exists $B'=(B'_{ij})_{i,j=1}^n\in  M_n(\mathcal A')$ such that
\begin{equation*}
A_{11}'\xi = (A'\tilde{\xi})_{1} = (B'\tilde{\eta})_{1}= \sum_{k=1}^n
B_{1k}'\eta_k \in \mathcal E,
\end{equation*}
so $A_{11}'\in \mathcal A'$ satisfies the theorem in the second case.

As the third and final case assume that $\mathcal E$ is an arbitrary subspace
of $H$ affiliated with $\mathcal A$. Define $\mathcal E_J = \sum_{\eta\in J}
\mathcal A'\eta$ for all finite subsets, $J$, of $\mathcal E$. Then
$(\mathcal E_J)$ is a upward filtering net ordered by inclusion with
upper bound $\mathcal E$, so if $p_J$ denotes the projection onto the
closure of the subspace $\mathcal E_J$ and $p$ denotes the projection
onto the closure of the subspace $\mathcal E$ then $(p_J)$ is an
increasing net of projections that converges strongly to $p$.

Let $\xi\in \overline{\mathcal E}$ and let $q$ be the projection onto
the closure of the subspace $\mathcal A'\xi$. Then $q\in \mathcal A$
and since $\overline{\mathcal A'\xi}\subseteq \overline{\mathcal E}$ we
have $q\leq p$ so $p_J\wedge q\nearrow q$ by lemma \ref{projectionlemma}. 
Choose $J$ a
finite subset of $\mathcal E$ such that 
\begin{equation} \label{delta1}
\norm{(p_J \wedge q)\xi -\xi}<\frac{\delta}{3}. 
\end{equation}
We know that $(p_J\wedge q)\xi \in \overline{\mathcal A'\xi}$ so the first part of the proof gives us an $A'\in \mathcal A'$ such that $A'(p_J\wedge q)\xi \in \mathcal A'\xi$ and 
\begin{equation}\label{delta2}
  \norm{A'(p_J\wedge q)\xi - (p_J\wedge q)\xi}<\frac{\delta}{3}.
\end{equation}
Choose $B'\in \mathcal A'$ such that $B'\xi = A'(p_J\wedge q)\xi$. Since 
\begin{equation*}
  B'\xi = A'(p_J\wedge q)\xi = (p_J\wedge q)A'\xi \subseteq p_J(H) = \overline{\sum_{\eta\in J} \mathcal A'\eta}
\end{equation*}
the second part of the proof gives us a $C'\in \mathcal A'$ such that
\begin{equation}\label{delta3}
  \norm{C'B'\xi - B'\xi}<\frac{\delta}{3}
\end{equation}
and $C'B'\xi\in \sum_{\eta\in J}\mathcal A'\eta\subseteq \mathcal E$. Combining (\ref{delta1}), (\ref{delta2}) and (\ref{delta3}) we have 
\begin{multline*}
  \norm{C'B'\xi - \xi} \leq    \norm{C'B'\xi - B'\xi} \\ + 
                       \norm{A'(p_J\wedge q)\xi - (p_J\wedge q)\xi}+
                         \norm{(p_J \wedge q)\xi -q\xi} \\
                       <  \frac{\delta}{3} + \frac{\delta}{3} +
                       \frac{\delta}{3} =\delta, 
\end{multline*}
so $C'B'\in \mathcal A'$ proves the theorem in the third case.
\end{proof}

\begin{proof}[Proof of theorem \ref{affsubspaces}]
Let $\mathcal E$ and $\mathcal F$ be subspaces of $H$ affiliated with
the finite von Neumann algebra $\mathcal A$ as stated in the
theorem. Let $\xi \in \overline{\mathcal E}\cap\overline{\mathcal
F}\subset \overline{\mathcal E}$. By theorem \ref{hjaelpetheorem} there
exists an $A'\in \mathcal A'$ such that $\norm{A'\xi-\xi}<
\frac{\delta}{2}$ for 
arbitrary $\delta>0$ and such that $A'\xi\in \mathcal E$ and $A'\xi \in
\overline{\mathcal E}\cap \overline{\mathcal F} \subseteq
\overline{\mathcal F}$. Applying theorem \ref{hjaelpetheorem} again we
obtain a $B'\in \mathcal A'$ such that
$\norm{B'(A'\xi)-A'\xi}<\frac{\delta}{2}$ and such that $B'(A'\xi)\in \mathcal F$. Since $\mathcal E$ is affiliated with
$\mathcal A$ we have $B'(A'\xi)\in \mathcal E$ and thus $B'A'\xi \in
\mathcal E\cap \mathcal F$. The inequalities imply that
\begin{equation*}
\norm{B'A'\xi - \xi} \leq \norm{B'A'\xi - A'\xi} +\norm{A'\xi-\xi} <
 \frac{\delta}{2} +  \frac{\delta}{2} = \delta.
\end{equation*}
Since $\delta>0$ is arbitrary we conclude that $\overline{\mathcal
E}\cap \overline{\mathcal F}\subseteq \overline{\mathcal E\cap\mathcal
F}$. Conversely $\mathcal E\cap \mathcal F \subseteq \overline{\mathcal
E}\cap \overline{\mathcal F}$ which is closed in $H$ so
$\overline{\mathcal E\cap \mathcal F} \subseteq \overline{\mathcal
E}\cap\overline{\mathcal F}$.
\end{proof}

\end{document}